  \documentclass[authoryear,preprint,3p,singlecolumn,11pt]{elsarticle}

\usepackage{xcolor}
\usepackage{amsthm}
\usepackage{subfigure}
\usepackage{graphicx}
\usepackage{caption}
\usepackage{amssymb}
\usepackage{lineno}
\usepackage{epsfig}
\usepackage[fleqn]{amsmath}
\usepackage{verbatim}
\newcommand{\eat}[1]{}

\newtheorem{Lemma}{Lemma}
\newtheorem{definition}{Definition}
\usepackage{tabulary}
\usepackage{booktabs}
\usepackage[ruled]{algorithm2e}
\usepackage{adjustbox}
\usepackage{rotating}
\usepackage{multirow}
\usepackage{lscape}
\usepackage{soul}
\usepackage{url}
 \usepackage{setspace}

 \newtheorem{Corollary}{Corollary}

\usepackage[noend]{algpseudocode}

\makeatletter
\def\ps@pprintTitle{
   \let\@oddhead\@empty
   \let\@evenhead\@empty
   \def\@oddfoot{\reset@font\hfil\thepage\hfil}
   \let\@evenfoot\@oddfoot
}
\makeatother

\journal{Transportation Research Part E}

\begin{document}
\begin{frontmatter}


\title{	Profit-Maximizing Parcel Locker Location Problem under Threshold Luce Model }



 \author[NUS]{Yun~Hui~Lin}
 \ead{isemlyh@gmail.com}
 
 \author[SUSS]{Yuan~Wang}
 \ead{wangyuan@suss.edu.sg}

 \author[NUS]{Loo Hay~Lee}
 \ead{iseleelh@nus.edu.sg}

 \author[NUS]{Ek Peng~Chew}
 \ead{isecep@nus.edu.sg}

 \address[NUS]{Department of Industrial Systems Engineering and Management, National University of Singapore}

 \address[SUSS]{School of Business, Singapore University of Social Sciences}

\begin{abstract}

The growth of e-commerce has created increasing complexity in logistics services. To remain competitive, logistics and e-commerce companies are exploring new modes as supplements to traditional home delivery, one of which is the self-service parcel locker.  This paper studies a parcel locker location problem where a company plans to introduce the locker service by locating locker facilities to attract customers. The objective is to maximize the profit, accounting for the revenue and the cost of facilities. To estimate the revenue, we use the threshold Luce model (TLM) to predict customers' likelihood of using the locker service. We then propose a combinatorial optimization model and develop exact solution methodologies that are practically implementable according to our extensive computational experiments. In effect, our modeling framework generalizes the traditional facility location problems based on the binomial logit model (BNL) and the multinomial logit model (MNL), both of which impose strong and strict assumptions on the customer's choice sets. That is, they assume that the choice sets will either contain only one facility or all facilities.  In our numerical experiment, we demonstrate that using the BNL and the MNL in the locker location problem could lead to, respectively, pessimistic and optimistic revenue estimation. Consequently, the suggested location decisions will be either conservative or  aggressive. Our proposed model, by contrast, can effectively relax these assumptions. Our results also reveal that the aggressive decision due to the use of the MNL will incur an unnecessarily high facility cost that cannot be compensated by the additional revenue, leading to profit loss that can be significant in various scenarios. Finally, we conduct sensitivity analysis on the input parameters and draw additional implications.




\end{abstract}
\begin{keyword}
Parcel locker \sep  Last-mile delivery  \sep  Facility location \sep Discrete choice model \sep Choice set restriction


\end{keyword}

\end{frontmatter}

\section{Introduction}
\label{intro}

E-commerce has been expanding at a remarkable rate as many customers turn to online shopping. The growth of e-commerce has led to a sharp increase in the delivery of online purchasing parcels.
Although logistics services have undergone significant improvements towards a more efficient and cost-effective manner, last-mile delivery (LMD) is still not efficient. According to Statista, a global business data platform, LMD accounts for 41\% of the overall supply chain cost in 2018\footnote{\tiny{https://www.statista.com/statistics/1043253/share-of-total-supply-chain-costs-by-type-worldwide/}}. The traditional approach for LMD is to deliver the parcel to the customer's doorstep. However, such an approach comes with various problems. Delivering packages to customers that are not present can cause delays on delivery routes or even lead to an empty trip. Such failed deliveries result in a waste of time and money for the company, and a hassle for the customer to reschedule the delivery that may create an unpleasant online shopping experience. Meanwhile, getting orders to individual addresses often involves one package per door, resulting in extra vehicle routes and higher fulfillment costs. To address the issues, logistics and e-commerce companies are constantly innovating in new delivery modes. A viable alternative is the usage of self-service lockers, which are akin to extra large mailboxes. With lockers, logistics companies can deliver the parcels in batch. The consolidation of shipments reduces vehicle movements and the number of vehicles required, offering companies better routes and improving cost efficiency.  The locker allows customers to select any locker location as a delivery address and then retrieve the parcels at their convenient time. Therefore, the locker service has been gaining popularity among logistics and e-commerce companies~\citep{Dieke2013,Morganti2014}, and it has been shown to improve the online shopping experience as well~\citep{Vakulenko2018}.

There is a proliferation of research on parcel lockers. Most of the existing works were empirical studies that focused on investigating the economic and environmental aspects of using parcel lockers~\citep{lachapelle2018parcel,Lemke2016,Punakivi2001,Song2013} or examining the customers' adoption and experience to this emerging technology~\citep{Vakulenko2018,yuen2018investigation,zhou2020understanding}. For example, \cite{Punakivi2001} conducted a case study in Finland and showed that the operational cost of home delivery is much more expensive than that of delivery to reception boxes and parcel lockers. \cite{lachapelle2018parcel} investigated the locker systems in 5 cities in Australia and discussed implications for urban and transport planning. 
\cite{Vakulenko2018} explored the experience of customer–parcel locker interaction in Sweden and customers' perceptions and value regarding parcel lockers, providing stakeholders with a foundation upon which to improve customers' parcel locker experiences.  Using survey data in Singapore, \cite{yuen2018investigation} argued that it is important is to integrate self-collection services into customers’ lifestyles, values, and needs to improve the adoption rate. A similar empirical study by \cite{zhou2020understanding} tested the influence of psychological factors and demographic information on behavioral intention to adopt the locker service.  

However, parcel locker facility location problems based on quantitative optimization approaches have received relatively limited attention.  Motivated by this, we study a locker location problem where a logistics company plans to introduce a parcel locker system by opening facilities to service customers. The objective is to maximize the profit, accounting for revenue and facility costs. In this paper, we apply the threshold Luce model (TLM), a discrete choice model recently proposed by~\cite{echenique2019general}, to predict the probability of a customer using the locker service. The revenue is then estimated using the predicted probability. In our setting, the revenue depends on the location of locker facilities.  An intuitive explanation is that if we open a locker near a customer zone, then this locker is attractive, and customers are more willing to use the locker service compared with the case where lockers are located far away. As a result, the company increases revenue but entails facility costs. In this paper, we propose a combinatorial optimization problem as a decision support framework for such a problem, allowing the company to make the location decision in the presence of customer choices and assess the trade-off between the revenue and the cost. \\

\noindent \textbf{Literature review.}  Before highlighting the contributions, we review related works and illustrate the positioning of this paper. To date, research on the facility location problem of the locker service  appears to be scarce. To the best of our knowledge,~\cite{deutsch2018parcel} developed the first quantitative approach to determine the optimal locker location. They assumed that customers maximize a deterministic utility and formulated an optimization problem to maximize the profit of the operator. Assuming the unlimited capacity of lockers, they showed that the problem is equivalent to the well-known uncapacitated facility location problem and can be solved by off-the-shelf optimization solvers.~\cite{hong2019routing} proposed a model to jointly determine the location of pick-up points and the delivery route. The objective is to minimize the weighted sum of the travel distance for the delivery vehicle and the customer, the operating costs for the pickup points, and the penalty for late delivery. More recently,~\cite{schwerdfeger2020optimizing} studied the dynamic deployment of mobile parcel lockers during the day to increase service accessibility, and~\cite{orenstein2019flexible} integrated the routing decision into the design of a flexible parcel delivery system. 
 
The above works explicitly or implicitly assumed that customers either make deterministic choices when deciding where to place their parcels or allow the operator to make choices for them. In reality, customers take the initiative in e-commerce and determine whether to use the locker service and which locker to select according to their preferences. These preferences are, however, unknown and hard to observe by the operator. It is thus argued that customers' choices should be interpreted probabilistically~\citep{gul2014random,mcfadden1981econometric,rusmevichientong2010dynamic}. This argument has been supported by empirical evidence in both locker-related research and facility location problems~\citep{drezner2002validating,drezner2006derived,gunawan2020parcel,suhara2021validating}. Therefore, it is meaningful to consider the probabilistic choice. To date, only two studies have integrated discrete choice models in parcel locker location problems, and both approaches were based on the multinomial logit model (MNL). Assuming that the total demand volume is fixed, \cite{lin2020last} investigated a locker network redesign problem to maximize the overall service level of a locker system. Inspired by the proposal of the Singapore Locker Alliance Network,~\cite{lyu2019last} developed a robust framework for a locker location problem under uncertainty. The objective is to maximize the utilization rate of the locker service and to serve the public residents. 

This paper is also related to the literature on facility location problems where choice models or choice rules are employed to forecast customer's choices. This research area is often called the competitive facility location (CFL), which investigates an ``entrant"  company, intending to enter a market with an incumbent or prospective competitor.  The decision problem is to locate new facilities to maximize the market share (or profit), given the competition from the competitor.  In general, there are two main research streams of the CFL according to the type of competition. The first one studies the ``static" competition and assumes that the location of competitor's facilities will not change when the company enters the market~\citep{lanvcinskas2017improving,gentile2018integer,lin2021generalized}), whereas the second stream considers the competitor's reaction and describes the problem from a game-theoretical perspective. We suggest readers to~\citet{mallozzi2017spatial} for a comprehensive discussion. In our setting, when making the location decision for the locker facilities,  the company has to consider the situation where customers may not use the lockers and instead seek logistics services from other providers (i.e., the outside option). This implies that the company is implicitly competing with the outside option to attract customers. Since the attraction of the outside option is fixed, our problem is indeed a static CFL model.

In the CFL literature, two classes of choice models (rules) have been widely applied to forecast the probability of a customer patronizing a facility, i.e.,  \textit{proportional choice rule} and \textit{partially binary choice rule}. Specifically, the proportional choice rule assumes that customers will consider all open facilities and that the probability of a customer patronizing a facility is proportional to the facility's attraction. The well-known and popular MNL and the gravity model belong to this class~\citep{drezner2018competitive,haase2014comparison,kuccukaydin2011discrete,ljubic2018outer}. The partially binary choice rule, by contrast, assumes that customers will consider two options, i.e., the facility with the highest attraction and the outside option. Customers then select an option from these two, proportional to the attraction~\citep{biesinger2016models,fernandez2017new,hakimi1990locations}. Such a choice rule resembles the binomial logit model (BNL) where the consideration set of the customer consists of only two alternatives.  

We point out that the assumptions of both choice rules (i.e., the MNL and the BNL) on customers' consideration sets are rather strong: they assume that customers will either consider all open facilities or only visit one facility. In facility location problems, these assumptions are not necessarily valid since, on one hand, a customer may not visit all open lockers, especially when the number of lockers is large (e.g., the number of locker stations has exceeded 200 in Singapore and is expected to reach around 1000 by the end of 2021\footnote{\tiny{https://www.straitstimes.com/tech/tech-news/govts-pick-parcel-delivery-locker-network-launched-with-over-200-lockers-deployed}}).
On the other hand, customers may not only consider one locker, especially in situations where multiple lockers appear to have nearly the same attraction to them. Therefore, it is reasonable to develop an approach that relaxes these assumptions.  \\

\noindent \textbf{Our contributions.} Positioning within the literature of parcel locker location problems and CFLs,  our contributions are as 
follows. (i) We investigate a locker location problem under the TLM 
and propose a combinatorial optimization model to maximize profit, accounting for revenue and the cost of facilities. We point out that both the BNL and the MNL  are special cases of the TLM; therefore, using the TLM naturally leads to a general and flexible modeling framework for the locker location problem, which also enables the company to experiment with the location decisions under the BNL and the MNL simply by tuning problem parameters; (ii) Besides the location decision, we also consider restricting customers' choice sets. Using an illustrative example, we show that when all open facilities are available to customers, revenue could be suboptimal, and the company can be better off by appropriately specifying the subset of open facilities that customers can use. To the best of our knowledge, the use of TLM and the restriction on customers' choice sets have not been investigated in the research of parcel locker location problems and CFLs thus far; (iii) To solve the proposed combinatorial optimization problem, we first formulate the problem as integer programs and then develop an exact solution approach based on conic reformulation. Our extensive computational experiments demonstrate that our methodology is capable of handling large-scale instances and thus practically implementable; (iv) Finally, we present numerical experiments to compare the solution structure under different choice models and conduct sensitivity analysis to draw insightful observations. \\

The outline for the paper is as follows.  Section~\ref{S:tlm} introduces the formal definition of threshold Luce model, followed by Section~\ref{S:pd} which presents a locker location problem under such a choice rule. We develop solution methodologies in Section~\ref{S:solution} and then conduct computational studies in Section~\ref{S:computational} to demonstrate the effectiveness of the proposed methods.  Section~\ref{S:numerical} presents numerical studies to compare the solutions under different choice models and conducts sensitivity analysis to investigate the impacts of problem parameters. Finally, we conclude the paper in Section~\ref{S:conclusion}.

\section{Threshold Luce model}
\label{S:tlm}

The threshold Luce model (TLM), proposed by~\cite{echenique2019general}, extends the classical proportional choice rule using the concept of \textit{dominance}. In its definition, if alternative $y$ dominates alternative $x$, denoted by $y  \succ x$, then the probability of selecting $x$ is zero in the presence of $y$. 

The TLM is essentially a special case of the general Luce model (GLM)~\citep{ahumada2018luce,echenique2019general}. Let $J$ denote the set of all alternatives and $a_j > 0$ be the attraction of alternative $j \in J$. In~\cite{flores2017assortment}, a non-select (outside) option, with index 0 and $a_0 > 0$, is introduced to capture the scenario where customers may not select any alternative. According to the GLM, given a subset of alternatives $S \subseteq J$,  a stochastic choice function $p$ returns a probability distribution over $c(S) \cup \{0\}$, where $c(S)$ is the set of all nondominated alternatives in $S$. Formally,  
\begin{definition}\citep{echenique2019general} 
\textit{A stochastic choice function $p$ is called a general Luce function if there exists a pair $(a,c)$ of attraction function $a \cup \{0\}: J \rightarrow R_{+}$ and a set function $c:2^{ J} \backslash \emptyset  \rightarrow 2^{J} \backslash \emptyset $ such that $c(S) \subseteq {S}$ for any subset $S$ of $J$, and the probability of selecting $j$ from set $S$ is given by}
\begin{align}
p_{j}(S)=\left\{\begin{array}{cc}{\frac{a_{j}}{\sum_{l \in c(S)} a_{l}+a_0}} & {\text { if } j \in c(S)} \\ {0} & {\text { if } j \notin c(S)}\end{array}\right. \quad \forall S \subseteq J.
\end{align}
\end{definition}
The function $c$ provides a way to capture the support of the stochastic choice function $p$. It prevents the dominated alternatives from being selected. As highlighted in~\cite{echenique2019general}, if $c(S) =S$, then $p$ coincides with the MNL. We further point out that if $|c(S)| =1$, then $p$ recovers the BNL.

In the GLM, the function $c$ does not necessarily depend on the attraction $a$. However, it is reasonable to assume that the dominance may be tied to attraction in some applications, 
leading to the following definition of the TLM.
\begin{definition} 
\textit{A stochastic choice function $p$ is called a threshold Luce function if there exists a nonnegative number $\gamma$ such that}
\begin{align}
c(S) =\left\{ j \in S \mid  (1+\gamma)a_{j} \geq a_{k}, \forall k \in S \right\}
\end{align}
\end{definition}
Here, the function $c$ is defined by $a$ and a nonnegative threshold parameter $\gamma$. The number $\gamma$ captures a threshold beyond which alternatives become dominated. That is, if $a_k > (1+\gamma)a_j$, then $j$ is dominated by $k$, and thus $j \notin c(S)$.  In other words, a ratio of more than $1+\gamma$ means that a less preferred alternative is dominated by a more preferred alternative, and the probability of selecting the dominated one will be zero. 

\section{Problem description}
\label{S:pd}

Having introduced the TLM, we proceed to the description of our locker location problem. Table~\ref{tab:Nomenclature} gives the main notations. Additional notations will be introduced when necessary. 

Consider a logistics company that plans to introduce a self-service locker system by opening parcel locker facilities, selected from some predetermined sites (denoted by set $J$), to service certain customer zones (denoted by set $I$). Each locker facility, once open, has an intrinsic attraction to customers. Meanwhile, customers may not select any locker to place their parcels and, instead, seek logistics services from other providers. We call this alternative ``the outside option". We assume that when the company sets up a chain of locker facilities and specifies the subset of open facilities that customers at a zone can use, customers determine whether to use the locker service from the company and which locker to patronize, following some discrete choice model. In our setting, the choices of customers are probabilistic. This is because customers are taking the initiative in e-commence, and they may prefer an option that is not necessarily with the highest utility according to their (unknown) preferences. Empirical results have supported that individual choice behaviors are typically probabilistic~\citep{gul2014random,gunawan2020parcel}. In this paper, we apply the TLM to predict the probability of a customer using the locker service and estimate the revenue. Then, the decision problem faced by the company is where to open the locker facilities and how to determine the subset of open lockers that customers can use so as to maximize the expected profit, accounting for the cost of locker facilities and the expected revenue under the TLM. 

\begin{table}[]

\caption{Notations.} \label{tab:Nomenclature}
\centering
\small
\begin{tabular}{lll}
\toprule
\noalign{\smallskip}
\multicolumn{3}{l}{\underline{Sets}}\\
$I$ &      $:$ & set of customer zones.\\
$J$ &     $:$ & set of candidate locker facilities.\\
$\{0\}$ &    $:$ & the outside option\\
$\Omega_{ij}$ & $:$ &set of lockers that are dominated by locker $j$ for customer zone $i$, $\forall i \in I, j \in J$. \\
\multicolumn{3}{l}{\underline{Parameters}}\\
$a_{ij}$ &     $:$ &attraction of option $j$ for customer zone $i$. $a_{ij} > 0$, $\forall j \in J \cup \{0\}$.\\
$d_i$ &     $:$ &demand at customer zone $i$, $\forall i \in I$. \\
$f_j$ &     $:$ & cost of locker facility $j$, $\forall j \in J$.\\
$\gamma$ &     $:$ &dominance threshold in threshold Luce model.\\
\multicolumn{3}{l}{\underline{Variables}}\\
$x_j$ & $:$ & 1, if locker $j$ is open; 0, otherwise, $\forall i \in I$.\\
$y_{ij}$ & $:$ & 1, if the open locker $j$ is in the nondominated set for customer zone $i$; 0, otherwise, $\forall i \in I, j \in J$.\\
         \bottomrule
\end{tabular}
\end{table}

Suppose the set of open lockers is denoted by $S$ ($S \subseteq J$) and the company entails a cost of locker facilities $F(S) = \sum_{j \in S} f_j$, where $f_j$ is the cost of facility $j$. We call $S$ the \textit{location decision}. We use $S_i$ to represent the set of lockers that the company allows customers at zone $i$ to select. Naturally,  $S_i \subseteq S \subseteq J$. The attraction of locker $j$ to customer zone $i$ is denoted by $a_{ij}, \forall i \in I, j \in J$. We require $a_{ij}$ to be strictly larger than 0. According to the TLM, the probability of a customer at zone $i$ selecting locker $j$ is
\begin{align}\label{eqt:p_tlm}
p_{ij}=\left\{\begin{array}{cc}{\frac{a_{ij}}{\sum_{l \in c_i(S_i)} a_{il}+ a_{i0}}} & {\text { if } j \in c_i(S_i)} \\ {0} & {\text { if } j \notin c_i(S_i)}\end{array}\right.,~\forall j \in S_i, i \in I
\end{align}
where $a_{i0}$ is the attraction of the outside option $0$, $\sum_{l \in c_i(S_i)} a_{il}$ stands for the attraction of the locker service, and the set $c_i(S_i)$ is defined as
\begin{align}
c_i(S_i) =\left\{ j \in S_i \mid  (1+\gamma)a_{ij} \geq a_{ik}, \forall k \in S_i \right\}
\end{align}
or equivalently,
\begin{align}
c_i(S_i) =\left\{ j \in S_i\mid \nexists k \in S_i:  (1+\gamma)a_{ij} < a_{ik} \right\}
\end{align}
The interpretation of the choice probability is that when customers are allowed to use a set of lockers $S_i$, they first discard all dominated lockers in $S_i$ and then select one from the nondominated set, $c(S_i)$, following the standard MNL with an additional outside option. Given $p_{ij}$, the expected revenue captured by the set of open lockers is
\begin{align}
R(S) =  \sum_{i \in I} \sum_{j \in S_i} d_i p_{ij}
\end{align}
where $d_i$ is the demand at customer zone $i, \forall i \in I$.  The problem is to find the subset $S$ of $J$ and to determine $S_i$ so that the profit $P(S)$ is maximized, i.e., 
\begin{align} \label{prob:abstract}
\max_{S_i \subseteq S\subseteq J}~&\{R(S) - F(S)\}
\end{align}
We refer to this combinatorial optimization problem as the \textit{locker location problem under the threshold Luce model} (LLPTL). 

In the above problem, the decision on $S_i$ reflects the ``restriction" on customers' choice sets, i.e., $S_i$ specifies the subset of open facilities that customers at zone $i$ can use, which is to be determined through optimization. We emphasize that  if we do not consider $S_i$ (i.e., set $S_i = S$), then opening an additional locker could lead to reduced revenue. This phenomenon reflects the \textit{paradox of choice} or \textit{choice overload effect}~\citep{chernev2015choice}.  
We present an example that illustrates such a scenario. \\

\noindent \textbf{Example.}  Consider two customer zones with demand $d_1 = d_2 = 50$, and two lockers with attraction $a_{11} = a_{12} = a_{21} = a_{22} = 2$. Assume the attraction of the outside option is $a_{10} = a_{20} = 4$ and $\gamma = 0.5$. Then $c_1(\{1,2\}) = \{1,2\}$ and  $c_2(\{1,2\}) = \{1,2\}$. In this case, the TLM is equivalent to the MNL and the choice probability is $p_{11} = p_{12} = p_{21} = p_{22} = 0.25$ and $p_{10} = p_{20} = 0.5$. Let $R_s$ and $R_t$ denote the revenue under the MNL and the TLM, then
\begin{align}
R_s = R_t = 50
\end{align}
Now, suppose a highly attractive locker, locker 3, is open with attraction $a_{13} = a_{23} = 3.1$.  By the MNL, $p_{11} = p_{12} = p_{21} = p_{22} = 0.18$, $p_{13} = p_{23} = 0.279$ and $p_{10} = p_{20} = 0.361$. Then,
\begin{align}
R_s^* = 63.9 > R_s
\end{align}
i.e., the revenue \textit{increases}. Meanwhile, $p_{10}$ and $p_{20}$ decrease. 

By the TLM, for customer zone 1, locker 3 dominates locker 1 and locker 2 because $3.1 > (1+\gamma) \times 2$. Therefore, $c_1(\{1,2,3\}) = 3$, meaning that $p_{11} = p_{12} = 0$, $p_{13} = 0.437$ and $p_{10} = 0.563$. Similarly, $c_2(\{1,2,3\}) = 3$, $p_{21} = p_{22} = 0$, $p_{23} = 0.437$ and $p_{20} = 0.563$. Then,
\begin{align}
R_t^*  = 43.7 < R_t
\end{align}
i.e., the revenue \textit{decreases}. $p_{10}$ and $p_{20}$ increase. The results are completely different from that by the MNL. Compared the values of $R^*_s$ and $R^*_t$, we find that the revenue estimated by the MNL is prominently higher. \hfill{$\square$}	 \\

The example shows that adding the highly attractive locker 3 can lead to reduced revenue under the TLM. However, if customers can only select either locker 1 or locker 2, then the revenue will not decrease. The implication is that, to maximize the revenue/profit, it is important to appropriately restrict the subset of open lockers that customers can use, instead of allowing customers to make the choices based on all open facilities. This observation justifies our proposed LLPTL model, which considers the decision on restricting customer's choice sets. \\

\noindent \textbf{Remark 1.} It is worth noting that when the company allows the customers to select lockers from all open facilities, the location problem should be modeled as a bilevel program or the so-called ``leader-follower-game" where the leader (the company) makes the location decision and the follower (the customer) acts accordingly based on the TLM. Obviously, the bilevel program will yield a lower profit since customer's choices are not restricted, and similar scenarios in the above example may occur. \\

\noindent \textbf{Remark 2.}  The choice set restriction is a somewhat new concept in facility location problems, and we do not notice any company that has practiced  it so far. Remaining on a conceptual level, a natural question then arises: ``\textit{How can the restriction be implemented?}" As modern e-commerce relies heavily on online platforms, a restriction mechanism can be designed as follows. After a customer purchases products online and decides to use parcel lockers as the delivery mode, he/she is required to input the address. Based on the address, the platform then suggests a set of lockers and requires the customer to select one of them to place the parcel. In this way, the company can restrict the choice set but also retain some degree of flexibility for the customer. \\

\noindent \textbf{Remark 3.}  There are different interpretations of the model regarding different values of $\gamma$. When $\gamma = 0$, customers will only consider the locker with the highest attraction that is accessible to them and the outside option.  This choice model is referred to as the \textit{partial binary rule}~\citep{suarez2004competitive,biesinger2016models} and resembles the BNL~\citep{gunawan2020parcel}. When $\gamma = \infty$, all open lockers have positive probabilities of being selected. The model recovers a location problem with the standard MNL because all open lockers will be accessible to customers to maximize the profit, which resembles the maximum capture facility location problem studied by~\cite{ljubic2018outer} or the competitive facility location by~\cite{mai2020multicut}. Since both special cases are NP-hard problems, the LLPTL problem itself is also NP-hard. \\

\section{Solution methodology}
\label{S:solution}

This section studies the solution methodology for the LLPTL.  Note that the combinatorial optimization problem presented in (\ref{prob:abstract}) is not explicit because the set $c_i(S_i)$ is only known implicitly. In this section, we first show that  (\ref{prob:abstract}) can be formulated as two equivalent integer programs (IPs). The first one expresses the dominance relation in the TLM using a large number of pairwise constraints, whereas the second one utilizes an aggregated formula to effectively model the dominance relation with significantly fewer constraints, but the continuous relaxation of the second formulation is weaker. Then, we visualize the dominance relation as a directed acyclic graph and propose path-based inequalities to improve the performance of the second formulation. Finally, we develop a unified optimization approach based on conic programming for both IPs. 

\subsection{Model formulation}

In  (\ref{prob:abstract}), set $S$ can be uniquely represented by a binary variable $x \in \{0, 1\}^{n}$ such that $j \in S$ if and only if locker $j$ is open, i.e., $x_j = 1$.  Define a binary variable $y \in \{0,1\}^{mn}$, such that $y_{ij} = 1$, if locker $j$ is in the nondominated set $c_i(S_i)$, i.e., $j \in c_i(S_i)$;  and $y_{ij} = 0$, if locker $j$ is dominated or not open, i.e., $j \notin c_i(S_i)$.  With $x_j$ and $y_{ij}$, we can write $F(S) = \sum_{j \in J} f_jx_j$ and express $R(S)$ as
\begin{align}
R(y) = \sum_{i\in I}d_i \frac{\sum_{j \in J}a_{ij} y_{ij}}{\sum_{j \in J}a_{ij} y_{ij} + a_{i0}} 
\end{align}
where $x_j$ and $y_{ij}$ should satisfy
\begin{align}
y_{ij} \leq x_j, \forall i \in I, \forall j \in J
\end{align}
which states that locker $j$ cannot be in the nondominated set if locker $j$ is not open. 

Now, we derive a simple inequality to represent the dominance relation between lockers.  Let $\Omega_{ij}$ be the set of lockers that are dominated by locker $j$ for customer zone $i$, i.e.,
\begin{align} \label{eqt:def_delta}
\Omega_{ij} = \{k \in J \mid a_{ij} > (1+\gamma)a_{ik} \}, \forall i \in I
\end{align}
Consider the following inequality 
\begin{align} \label{eqt:dominated}
y_{ij}+y_{ik} \leq 1,  \forall k \in \Omega_{ij}, \forall i \in I, j \in J
\end{align}
which enforces that at most one of $y_{ij}$ and $y_{ik}$ can take the value of 1 if locker $j$ dominates locker $k$. Since $a_{ij} > a_{ik}$, this inequality effectively imposes that when both lockers are available to customers at zone $i$,  we have $y_{ik} = 0$ in the optimal solution. We refer to this constraint as  the \textit{disaggregated dominance constraint} (DDC) because the dominance relation is imposed in pairs. A similar constraint applied to the assortment optimization can be found in~\cite{flores2017assortment}.  With the DDC, we have the following IP to model the LLPTL in~(\ref{prob:abstract}). 
\begin{subequations}\label{model:MIP-d}
\begin{alignat}{4}
  \max_{x,y}~&\sum_{i\in I}d_i \frac{\sum_{j \in J}a_{ij} y_{ij}}{\sum_{j \in J}a_{ij} y_{ij} + a_{i0}} - \sum_{j \in J} f_j x_j \\
\textbf{[IP-D]} \qquad \text{st}.~& y_{ij} \leq x_j,  \forall i \in I, j \in J \\
&  y_{ij}+y_{ik} \leq 1,   \forall k \in \Omega_{ij}, \forall i \in I, j \in J \\
& x_j \in \{0,1\}, \forall j \in J \\
& y_{ij} \in \{0,1\}, \forall i \in I, j \in J
\end{alignat}
\end{subequations}

\eat{
The following results show the hardness of the problem.
\begin{Lemma} \label{prop:np-hard}
MIP-d is NP-hard.
\end{Lemma}

\begin{proof}
To show the NP-hardness, we show that MIP-d is NP-hard even when there is only 1 customer zone.   Suppose $|I| =1$, then $d=1$ and we can drop the subscript $i$ in MIP-d. Since we will open $P$ facilities, then by $y_j \leq x_j$, we have $\sum_{j \in J} y_j \leq P$. Meanwhile, the objective function $\frac{\sum_{j \in J}a_{j} y_{j}}{\sum_{j \in J}a_{j} y_{j} + a_{0}}$ is an increasing function due to $a_j \geq 0, \forall j \in J$. To maximize this function is equivalent to maximize $\sum_{j \in J}a_j y_j$. Therefore, MIP-d with $|I|=1$ can be represented by
\begin{equation} \label{model:MIP-d_reduced_eq}
\begin{aligned} 
\max_{y}~& \sum_{j \in J}a_{j} y_{j} \\
\text{st}.~& \sum_{j \in J} y_j \leq P \\
&  y_{j}+y_{k} \leq 1, \forall j \succ k \\
& y_j \in \{0,1\}, \forall j \in J
 \end{aligned}
\end{equation}
which can be visualized as the maximum weighted independent set problem of size at most $P$ for a bipartite graph (MWISP), i.e., MWISP for a graph $G = (V=V_1 \cup V_2, E)$, where $V_1 \cap V_2 =  \emptyset$, every edge $(j,k) \in E$ satisfies $j \in V_1$ and $k \in V_2$, $u_j$ is the weight of vertex $j$ and $P$ is the budget. It has been proven that MWISP is NP-hard~\citep{kalra2017maximum}; therefore, MIP-d is NP-hard as well.
\end{proof}

In fact, even when the location decision has been made, the remaining problem of [MIP-d] with only $y$ variable is still NP-hard. This argument is justified in the following lemma.
\begin{Lemma}
Given a feasible location decision $\bar{x}$, determining the optimal assignment decision $y$ is a NP-hard problem.
\end{Lemma}
\begin{proof}
Define set $\Xi$ as the set of open lockers suggested by $\bar{x}$, i.e., $\Xi = \{j \mid \bar{x}_j = 1, \forall j \in J\}$. Similar to the proof in Lemma~\ref{prop:np-hard}, we consider the reduction of the problem to the case where $|I| = 1$. Then the optimal assignment decision can be obtained by
\begin{equation} \label{model:MIP-d_reduced_eq}
\begin{aligned} 
\max_{y}~& \sum_{j \in \Xi}a_{j} y_{j} \\
\text{st}.~&  y_{j}+y_{k} \leq 1, \forall j \succ k \\
& y_j \in \{0,1\}, \forall j \in \Xi
 \end{aligned}
\end{equation}
which is equivalent to the well-known NP-hard problem of finding the maximum independent set for a  graph with $|\Xi|= P$ vertices.
\end{proof}
}

In the DDC, the number of constraints can go up to $O(mn^2)$ because the locker with the highest attraction can dominate at most $n-1$ lockers; the locker with the second-highest attraction can dominate at most $n-2$ locker; and so on. The number of constraints can be extremely large in some scenarios. In light of this, we present an alternative representation of the dominance with $O(mn)$ constraints. 

Observe that if $y_{ij} = 1$ (i.e., locker $j$ is open and in the nondenominated set for customer zone $i$),  then those lockers that are dominated by $j$ will not be in the nondominated set (i.e., $\sum_{k \in \Omega_{ij}} y_{ik} = 0$). We can express this condition as
\begin{align}
\label{eqt:dominated_a}
\sum_{k \in \Omega_{ij}} y_{ik} \leq |\Omega_{ij}| \cdot (1-y_{ij}) , \forall i \in I, j \in J
\end{align}
We refer to this constraint as the \textit{aggregated dominance constraint} (ADC) to reflect the fact that we aggregate the pairwise dominance relations of locker $j$ for customer zone $i$ into a single constraint. The following lemma establishes the relationship between DDC and ADC. 
\begin{Lemma} \label{prop:DDC_ADC}
DDC yields a stronger continuous relaxation than ADC.
\end{Lemma}

\begin{proof}
 Write (\ref{eqt:dominated}) as
\begin{align}
y_{ik} \leq 1 - y_{ij},  \forall  k \in  \Omega_{ij}, \forall i \in I, j \in J
\end{align}
Take the summation on both sides with respect to $k$, and we have:
\begin{align}
\sum_{k \in \Omega_{ij}}y_{ik} \leq \sum_{k \in \Omega_{ij}}(1 - y_{ij}) = |\Omega_{ij}|\cdot (1-y_{ij}), \forall i \in I, j \in J
\end{align}
which is exactly the ADC  in (\ref{eqt:dominated_a}). This indicates the polyhedron defined by the continuous relaxation of (\ref{eqt:dominated}) is a subset of that by (\ref{eqt:dominated_a}). Therefore, the DDC can yield a tighter relaxation bound than the ADC~\citep{chen2011applied}. 
\end{proof}

Although replacing the DDC by the ADC can significantly reduce the number of constraints, the resultant continuous relaxation will become weaker as suggested by the above lemma. 

\subsection{Improving the ADC with path-based inequalities}

\begin{figure*}[h]
\begin{center}
     \psfig{figure=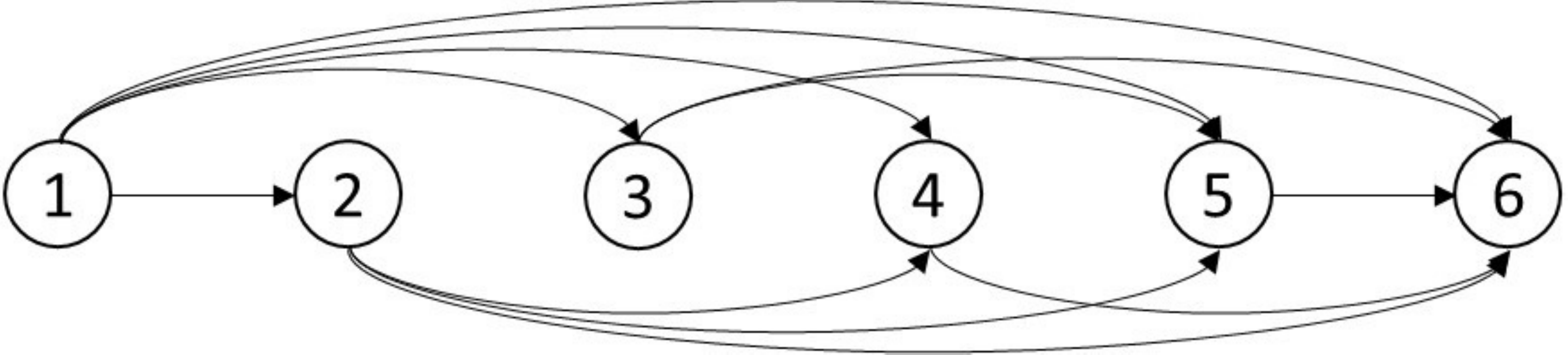,width=100mm}
\end{center}
\caption{A directed acyclic graph with 6 vertices.} \label{fig:DAG}
\end{figure*}

Due to weak relaxation of the ADC, we seek methods to strengthen the formulation and to expedite the computation. The dominance relation is a partial order that is irreflexive, transitive and antisymmetric~\citep{flores2017assortment}. For customer zone $i$, the dominance relation can be represented as a Directed Acyclic Graph (DAG), in which vertices represent the lockers, and there is a directed edge $(m,n)$ if and only if $m \succ n$. Formally, let $G_i(J,E_i)$ be the DAG associated with customer zone $i$, where $J$ is the vertex set (i.e., the set of candidate lockers) and $E_i$ is the edge set. An edge $(j,k)$ exists in $G_i$ if and only if $a_{ij} >  (1+\gamma) a_{ik}$. For example, suppose that, for a customer zone, there are 6 lockers whose attraction is ranked from the largest to the smallest. The DAG, as in Figure~\ref{fig:DAG}, represents the following dominance relations: $1\succ \{2,3,4,5,6\}$; $2\succ \{4,5,6\}$; $3\succ \{5,6\}$; $4\succ 6$; $5\succ 6$.

A path in a DAG is a sequence of vertices that are connected by directed edges. In our problem, a path represents a sequence of dominance relations. For instance, a path $1\to5\to6$ in Figure~\ref{fig:DAG} with the vertices sequence $\{ 1, 5, 6\}$ indicates that $1 \succ 5 \succ 6$. With the path relation, when locker $1$ is open, locker $5$ and locker $6$ will be in the dominated set and only locker $1$ can have a positive probability of being selected. Similarly, when locker $1$ is not open and locker $5$ is open, locker $6$ will be in the dominated set and only locker $5$ can have a positive probability of being selected. Therefore, at most one of the lockers in the vertex sequence $\{ 1, 5, 6\}$ can have a positive probability. This interpretation implies that the dominance constraint can be imposed in terms of paths. Formally,
\begin{Lemma} \label{prop:valid}
Suppose there exists a path $r_i$ in $G_i$ that transverses vertex sequence $J_{r_i}$, where $J_{r_i} \subseteq J$, then
\begin{align} \label{eqt:path}
\sum_{j \in J_{r_i}} y_{ij} \leq 1
\end{align}
is valid.
\end{Lemma}
\begin{proof}
Suppose the vertex sequence is $\{k^1, k^2, k^3,...,k^t\}$ with $t = |J_{r_i}|$, then $k^1 \succ k^2 \succ k^3 \succ ... \succ k^t$. Assume locker $k^1$ is open, then $\sum_{j \in J_{r_i} \backslash	\{k^1\}} y_{ij} = 0$ since $k^1$ dominates all lockers in the set $J_{r_i} \backslash \{k^1\}$.  Now, assume locker $k^1$ is not open and $k^2$ is open, then $\sum_{j \in J_{r_i} \backslash \{k^2\}} y_{ij} = 0$ since $x_{k^1} = 0$ and $k^2$ dominates all lockers in the set $J_{r_i} \backslash	\{k^1,k^2\}$.  Similarly, assume locker $k^1$ and $k^2$ are not open and $k^3$ is open, then $\sum_{j \in J_{r_i} \backslash \{k^3\}} y_{ij} = 0$ since $x_{k^1} = x_{k^2} = 0$ and $k^3$ dominates all lockers in the set $J_{r_i} \backslash \{k^1,k^2,k^3\}$. Continue this procedure. We can conclude that $\sum_{j \in J_{r_i}} y_{ij} \leq 1$.
\end{proof}
We refer to (\ref{eqt:path}) as a \textit{path-based inequality}. It indicates that customers in zone $i$ will visit at most one locker in $J_{r_i}$. The following lemma establishes the linkage between such an inequality and the pairwise inequality in (\ref{eqt:dominated}). 

\begin{Corollary} \label{prop:path}
A path-based inequality in (\ref{eqt:path}) simutanenously implies ${t(t-1)}/{2}$ pairwise dominance relations, where $t = |J_{r_i}|$.
\end{Corollary}
\begin{proof}
This immediately follows from the observation that the first vertex dominates $t-1$ downstream vertices; the second vertex dominates $t-2$ downstream vertices; and so on.
\end{proof}

\begin{definition}
A path $r$ transversing vertex set $J_{r_i}$ in $G_i$ is called a \textit{maximal path} if there does not exist a distinct path $\hat{r}$ transversing vertex set $J_{\hat{r}_i}$ such that $J_{r_i} \subset J_{\hat{r}_i}$.
\end{definition}
Informally, a path is maximal if we cannot further add any connected edge from the DAG to it. The dominance relation embedded in $G_i $ can be fully represented by the set of all maximal paths because such a set will contain the information about all edges in $G_i$ and the correct pairwise relation. In principle, we can find out all maximal paths and impose the corresponding constraints. However, there is an issue: The number of edges in all $G_i$ is equal to the number of constraints in the DDC. There could be an exponential number of maximal paths. Exhaustive enumeration of them will incur significant computational time. 

To reduce the effort, we use one maximal path of each $G_i$. Specifically, we are interested in the path with the largest number of vertices, namely, the longest path. It can be derived using a linear-time algorithm starting from the \textit{topological sorting}, which sorts a DAG into a linear ordering of its vertices such that for every directed edge $(x,y)$ from vertex $x$ to vertex $y$, $x$ comes before $y$ in the ordering (e.g., in Figure~\ref{fig:DAG} where every edge points to the right side). We can then process all vertices one by one in topological order.  For each vertex $x$, compute the length of the longest path to $x$ by its predecessors and add one to the maximum length recorded for those predecessor. If $x$ has no predecessor, the length of the longest path to $x$ is zero. The longest path can then be obtained by starting at the vertex $x$ with the largest recorded value, then repeatedly tracing backwards to its predecessor with the largest recorded value, and reversing the sequence of vertices. Such an algorithm has been successfully implemented by various advanced computer graph and network packages. In this paper, we leverage the graph algorithm in the Python package NetworkX~\citep{schult2008exploring} so that we avoid the need for designing the algorithm from scratch.

Let $\bar{J}_i$ be the set of vertices in the longest path for $G_i$, then we impose the following path-based constraint for customer zone $i$
\begin{align} \label{eqt:long_path}
\sum_{j \in \bar{J}_i} y_{ij} \leq 1, \forall i \in I
\end{align}
Together with the ADC, we propose the following IP to model the LLPTL.
\begin{subequations}\label{model:MIP-s}
\begin{alignat}{4}
\max_{x,y}~&\sum_{i\in I}d_i \frac{\sum_{j \in J}a_{ij} y_{ij}}{\sum_{j \in J}a_{ij} y_{ij} + a_{i0}} - \sum_{j \in J} f_j x_j \\
\text{st}.~& y_{ij} \leq x_j,  \forall i \in I, j \in J \\
\textbf{[IP-A]}~\qquad & \sum_{k \in \Omega_{ij}} y_{ik} \leq  |\Omega_{ij}| \cdot (1-y_{ij}) , \forall i \in I, j \in J \\
&\sum_{j \in \bar{J}_i} y_{ij} \leq 1, \forall i \in I \\
& x_j \in \{0,1\}, \forall j \in J \\
& y_{ij} \in \{0,1\}, \forall i \in I, j \in J
\end{alignat}
\end{subequations}

\subsection{Solution approach}
\label{S:micqp}

In this subsection, we develop an exact approach to solve IP-D and IP-A. We first define the set $\Xi_A$ and $\Xi_D$ as the feasible region of IP-A and IP-D. We then use $\Xi$ to refer to $\Xi_A$ and $\Xi_D$. It is easy to see that the IP model is equivalent to the following  minimization problem
\begin{equation}\label{model:MIP-min}
\begin{aligned} 
\min_{x,y}~& \sum_{i\in I} \frac{d_i}{\sum_{j \in J}\pi_{ij} y_{ij} + 1} + \sum_{j \in J} f_j x_j \\
st.~& (x,y) \in \Xi
\end{aligned}
\end{equation}
where $\pi_{ij} = a_{ij}/a_{i0}$, and the objective reflects the summation of lost demand and the cost of facilities. Note that the objective function is convex since the first term is the composition of a convex and decreasing function, $1/z$ for $z>0$, with an affine function $ \sum_{j \in J}\pi_{ij} y_{ij}+ 1$. Motivated by the convexity, we present an approach to recast (\ref{model:MIP-min}) into a mixed-integer conic quadratic program (MICQP).

To proceed, let $z_i = \sum_{j \in J}\pi_{ij} y_{ij}+ 1$. The objective function of~(\ref{model:MIP-min}) becomes 
\begin{align}
\sum_{i\in I}d_i/z_i +  \sum_{j \in J} f_j x_j
\end{align}
Define $\beta_i$ such that $\beta_i \geq 1/z_i$.  We have $\beta_i z_i \geq 1, \forall i \in I $.~(\ref{model:MIP-min}) can then be restated as the following program:
\begin{subequations}\label{model:MICQP}
\begin{alignat}{4}
 \min_{x,y}~& \sum_{i \in I} d_{i} \beta_{i} + \sum_{j \in J} f_j x_j\\
\label{rotated_cone} \text {st.}~& \beta_{i} z_{i} \geq 1, \forall i \in I \\ 
\textbf{[MICQP]} \qquad & z_i = \sum_{j \in J}\pi_{ij} y_{ij}+ 1, \forall i \in I \\
&0 \leq \beta_{i} \leq 1, \forall i \in I \\
& (x,y) \in \Xi 
\end{alignat}
\end{subequations}
where the objective function and constraints are linear except constraint~(\ref{rotated_cone}). In fact, since $z_i \geq 0$ and $\beta_i \geq 0$, (\ref{rotated_cone}) is equivalent to
\begin{equation} \label{eqt:cone2}
|| (2,z_i-\beta_i)||_2 \leq z_i + \beta_i
\end{equation}
which is a conic quadratic inequality. Therefore, (\ref{model:MICQP}) is indeed a MICQP, which can be handled by advanced mixed-integer second-order conic solvers (e.g., CPLEX and Gurobi). It can also be verified that $\beta_iz_i = 1$  in the optimal solution. We can then enforce that $\beta_i \leq 1$ since  $z_i \geq 1$ by definition.  

\begin{Lemma} \label{prop:equivalent}
Problem (\ref{model:MIP-min}) is equivalent to problem (\ref{model:MICQP}).
\end{Lemma}
\begin{proof}
Let $(x^1,y^1)$ be the optimal solution to problem (\ref{model:MIP-min}) and $(x^2,y^2,z^*,\beta^*)$ be the optimal solution to problem(\ref{model:MICQP}). Denote the corresponding objectives as $r_1(\cdot)$ and $r_2(\cdot)$. We show that $r_1(x^1,y^1)= r_2(x^2,y^2,z^*,\beta^*)$.

Given the optimal solution $(x^1,y^1)$ to  problem (\ref{model:MIP-min}), a feasible solution to problem (\ref{model:MICQP}) can be constructed as $(x^1,y^1, \bar{z},\bar{\beta})$, where $\bar{z}_i = \sum_{j\in J} \pi_{ij} y^1_{ij} + 1$ and $\bar{\beta}_i = 1/\bar{z}_i $, $\forall i \in I$. Therefore, we have
\begin{align} \label{eqt:lhs}
r_1(x^1,y^1) =  r_2(x^1,y^1, \bar{z},\bar{\beta}) \geq r_2(x^2,y^2,z^*,\beta^*)
\end{align}
The equality follows from the construction of $(x^1,y^1, \bar{z},\bar{\beta})$ and the inequality comes from the optimality of $(x^2,y^2,z^*,\beta^*)$ to problem (\ref{model:MICQP}).

On the other hand, given the optimal solution $(x^2,y^2,z^*,\beta^*)$ to problem (\ref{model:MICQP}), we have $z^*_i = \sum_{j\in J} \pi_{ij} y^2_{ij} + 1$, $\forall i \in I$, then $r_2(x^2,y^2,z^*,\beta^*) = r_1(x^2,y^2)$. Therefore, we have
\begin{align}\label{eqt:rhs}
r_2(x^2,y^2,z^*,\beta^*) = r_1(x^2,y^2) \geq r_1(x^1,y^1)
\end{align}
The inequality is due to the optimality of $(x^1,y^1)$ to  problem (\ref{model:MIP-min}).

Combining (\ref{eqt:lhs}) and  (\ref{eqt:rhs}), we can conclude that $r_1(x^1,y^1)= r_2(x^2,y^2,z^*,\beta^*)$.
\end{proof}


\section{Computational experiment}
\label{S:computational}

This section presents the computational experiment and analysis. We refer to the MICQP approaches applied to IP-D and  IP-A as {\fontfamily{qcr}\selectfont{MICQP-D}} and {\fontfamily{qcr}\selectfont{MICQP-A}}, respectively. To illustrate the effectiveness of both approaches, we compare them with the mixed-integer convex modeling package CVXPY.  Using CVXPY, we directly model and solve IP-D. We refer to the use of  CVXPY as {\fontfamily{qcr}\selectfont{CVX-D}}. The computational experiments are done using the Python-Gurobi interface (version 9.1.2) on a 16 GB Memory Mac computer with a 2.6 GHz Intel Core i7 processor. Without loss of generality, we assume that the costs of locker facilities are the same, i.e., $f = f_j, \forall j \in J$.

\subsection{Experiment set up}

We use two artificially datasets to generate benchmark instances. They are described as follows\footnote{All data are generated with fixed random seeds and can be provided upon request.}.

\begin{itemize}
\item[\textbf{DS1}] A medium-scale dataset consists of 200 customer zones and 100 candidate locker facilities. The locations of the demand zones and the candidate facilities are generated from a 2-dimensional uniform distribution $[0,30]^2$.  The demand $d_i$ is generated from a uniform distribution $[1,1000]$.
\item[\textbf{DS2}] A large-scale dataset consists of 400 customers zones and 150 candidate locker facilities. The locations of the demand zones and the candidate facilities are generated from a 2-dimensional uniform distribution $[0,40]^2$. The demand $d_i$ is generated from a uniform distribution $[1,1000]$.
\end{itemize}
For both datasets, the distance matrix $L$ is measured by the Euclidean distance, and the attraction of locker $j$ to customer zone $i$ is computed by 
\begin{align} \label{eqt:attr_locker}
a_{ij} = e^{-\alpha L_{ij}}, \forall i \in I, j \in J
\end{align}
where $u_{ij} = -\alpha L_{ij}$ stands for the utility. Indeed, (\ref{eqt:attr_locker}) represents an exponential decay function that has been applied to various location problems~\citep{hodgson1981location,aboolian2007competitive,drezner2008lost}. It reflects the decline in the attraction of a facility to a customer as an exponential function of the distance between them. The parameter $\alpha$ controls customers' sensitivity to the distance. With an extremely large $\alpha$ value, customers will almost always patronize the nearest facility if they decide to use the locker service. In this section, $\alpha$ is fixed at 1 for parameter simplicity, but we will explore its impact on the solution structure by varying its value at next section. 

 For the outside option 0, we assume that the attraction is
\begin{align}\label{eqt:attr_outside}
a_{i0} = \xi \cdot e^{-1}, \forall i \in I
\end{align}
where $\xi >0$ is a parameter that controls the magnitude of the attraction. A larger value of $\xi$ indicates that the outside option is more attractive, and customers are more willing to select it. 

Finally, we introduce the notations that will be used in our later discussion.
\begin{itemize}
\item $t[s]$ denotes the computational time in seconds. An instance is labeled with ``t.l." when it cannot be solved by an approach within some predetermined time limit.
\item  $gp[\%]$ stands for the relative exit gap in percentages when the algorithm terminates. It is computed as  $|zbb-zopt|/|zbb| \times 100$, where $zopt$ is the value of the optimal integer solution and $zbb$ is the best bound at termination. If an instance is solved to optimality, then zbb equals zopt (or within the default optimality gap tolerance 0.01\%). 

\item  $\#N$ is the number of branch-and-cut nodes explored by Gurobi.
\end{itemize}

\subsection{Result analysis on DS1}

We first present the computational results on DS1 instances. For this dataset, the time limit is set to 3600 seconds. Table~\ref{Tab:Result200-1}-\ref{Tab:Result200-05} show the results under $\xi = 1$ and $\xi=0.5$ for different combinations of $(f,\gamma)$. For each problem instance, the computational time for the approach that performs the best (i.e., the smallest CPU time among {\fontfamily{qcr}\selectfont{CVX-D}}, {\fontfamily{qcr}\selectfont{MICQP-D}}, and {\fontfamily{qcr}\selectfont{MICQP-A}} to reach optimality) is highlighted in boldface. 

\begin{table}[h]
\footnotesize
\centering
\caption{Computational results of the DS1 instances under $\xi=1$.  Time limit is 3600 seconds.} \label{Tab:Result200-1}
\begin{tabular}{lrrrrrrrrrrrr} 
\toprule
\multirow{2}{*}{$f$} & \multirow{2}{*}{$\gamma$} & \multicolumn{3}{c}{{\fontfamily{qcr}\selectfont{CVX-D}}}                                                     &  & \multicolumn{3}{c}{{\fontfamily{qcr}\selectfont{MICQP-D}}}                                                     &  & \multicolumn{3}{c}{{\fontfamily{qcr}\selectfont{MICQP-A}}}                                                      \\ 
\cline{3-5}\cline{7-9}\cline{11-13}
                     &                           & \multicolumn{1}{c}{t[s]} & \multicolumn{1}{c}{gp[\%]} & \multicolumn{1}{c}{\#N} &  & \multicolumn{1}{c}{t[s]} & \multicolumn{1}{c}{gp[\%]} & \multicolumn{1}{c}{\#N} &  & \multicolumn{1}{c}{t[s]} & \multicolumn{1}{c}{gp[\%]} & \multicolumn{1}{c}{\#N}  \\ 
\hline
500     & 1             & 865.8                    & 0                          & 1864                    &  & \textbf{85.0}            & 0                          & 216                     &  & 189.7                    & 0                          & 2330                     \\
     & 2             & 1572.8                   & 0                          & 1633                    &  & \textbf{181.8}           & 0                          & 453                     &  & 317.8                    & 0                          & 3884                     \\
     & 3             & 1840.1                   & 0                          & 2405                    &  & \textbf{207.1}           & 0                          & 641                     &  & 427.6                    & 0                          & 4019                     \\
     & 5             & 1749.1                   & 0                          & 3732                    &  & \textbf{358.9}           & 0                          & 1114                    &  & 429.6                    & 0                          & 9379                     \\
     & 7             & 2211.1                   & 0                          & 9450                    &  & 813.1                    & 0                          & 690                     &  & \textbf{301.2}           & 0                          & 2489                     \\
     & 10            & 2010.3                   & 0                          & 8554                    &  & 367.5                    & 0                          & 1709                    &  & \textbf{252.7}           & 0                          & 2183                     \\
     & 15            & 1699.7                   & 0                          & 4791                    &  & 554.5                    & 0                          & 1045                    &  & \textbf{241.8}           & 0                          & 8853                     \\
     & 20            & 1503.6                   & 0                          & 5793                    &  & \textbf{362.1}           & 0                          & 1065                    &  & 490.7                    & 0                          & 10124                    \\
     & Avg  &    1681.6                       &        0                    &      4778                   &  &           366.3               &          0                  &             867            &  &                           \textbf{331.4} &                0            &           4490               \\
     &               &                          &                            &                         &  &                          &                            &                         &  &                          &                            &                          \\
1000     & 1             & 1469.2                   & 0                          & 905                     &  & 304.9                    & 0                          & 1529                    &  & \textbf{227.4}           & 0                          & 2889                     \\
     & 2             & 1503.3                   & 0                          & 1072                    &  & \textbf{352.9}           & 0                          & 2154                    &  & 611.0                    & 0                          & 7119                     \\
     & 3             & 1919.5                   & 0                          & 1115                    &  & 841.8                    & 0                          & 1503                    &  & \textbf{501.9}           & 0                          & 2736                     \\
     & 5             & 2548.2                   & 0                          & 1130                    &  & 795.1                    & 0                          & 1773                    &  & \textbf{626.2}           & 0                          & 3229                     \\
     & 7             & 3089.6                   & 0                          & 1439                    &  & 712.7                    & 0                          & 1404                    &  & \textbf{564.1}           & 0                          & 3304                     \\
     & 10            & 2331.9                   & 0                          & 1360                    &  & 633.4                    & 0                          & 941                     &  & \textbf{472.0}           & 0                          & 10052                    \\
     & 15            & 2066.7                   & 0                          & 1740                    &  & 681.5                    & 0                          & 958                     &  & \textbf{389.7}           & 0                          & 2357                     \\
     & 20            & 2155.8                   & 0                          & 1660                    &  & 677.1                    & 0                          & 1662                    &  & \textbf{404.2}           & 0                          & 3265                     \\
     & Avg &         2135.5                 &        0                    &     1303                    &  &    624.9                      &          0                  &     1491                    &  &    \textbf{474.6}                      &      0                      &     4369                     \\
\bottomrule
\end{tabular}
\end{table}

\begin{table}[h]
\footnotesize
\centering
\caption{Computational results of the DS1 instances under $\xi=0.5$. Time limit is 3600 seconds.}
\label{Tab:Result200-05}
\begin{tabular}{lrrrrrrrrrrrr} 
\toprule
\multirow{2}{*}{$f$} & \multirow{2}{*}{$\gamma$} & \multicolumn{3}{c}{{\fontfamily{qcr}\selectfont{CVX-D}}}                                                     &  & \multicolumn{3}{c}{{\fontfamily{qcr}\selectfont{MICQP-D}}}                                                     &  & \multicolumn{3}{c}{{\fontfamily{qcr}\selectfont{MICQP-A}}}                                                      \\ 
\cline{3-5}\cline{7-9}\cline{11-13}
                     &                           & \multicolumn{1}{c}{t[s]} & \multicolumn{1}{c}{gp[\%]} & \multicolumn{1}{c}{\#N} &  & \multicolumn{1}{c}{t[s]} & \multicolumn{1}{c}{gp[\%]} & \multicolumn{1}{c}{\#N} &  & \multicolumn{1}{c}{t[s]} & \multicolumn{1}{c}{gp[\%]} & \multicolumn{1}{c}{\#N}  \\ 
\hline
500                  & 1                         & 2327.4                   & 0                          & 6417                    &  & \textbf{120.6}           & 0                          & 323                     &  & 432.3                    & 0                          & 1354                     \\
                     & 2                         & t.l.                     & 0.32                       & 3866                    &  & \textbf{613.0}           & 0                          & 1591                    &  & 1132.5                   & 0                          & 5881                     \\
                     & 3                         &  t.l.                      & 0.50                       & 3354                    &  & \textbf{511.2}           & 0                          & 1424                    &  & 542.3                    & 0                          & 2735                     \\
                     & 5                         &  t.l.                      & 0.20                       & 7841                    &  & 740.3                    & 0                          & 736                     &  & \textbf{501.9}           & 0                          & 1866                     \\
                     & 7                         &  t.l.                     & 0.18                       & 9137                    &  & 705.3                    & 0                          & 1417                    &  & \textbf{506.7}           & 0                          & 3368                     \\
                     & 10                        &  t.l.                      & 0.16                       & 16229                   &  & 659.3                    & 0                          & 1657                    &  & \textbf{344.5}           & 0                          & 3135                     \\
                     & 15                        &  t.l.                      & 0.10                       & 21924                   &  & 537.2                    & 0                          & 1035                    &  & \textbf{378.5}           & 0                          & 2214                     \\
                     & 20                        &  t.l.                      & 0.06                       & 21493                   &  & 468.5                    & 0                          & 936                     &  & \textbf{242.2}           & 0                          & 6256                     \\
                     &                           &                          &                            &                         &  &                          &                            &                         &  &                          &                            &                          \\
1000                 & 1                         &  t.l.                      & 0.80                       & 1857                    &  & \textbf{889.5}           & 0                          & 1397                    &  & 1105.3                   & 0                          & 3162                     \\
                     & 2                         &  t.l.                      & 1.19                       & 634                     &  & \textbf{1815.1}          & 0                          & 3363                    &  & 2808.5                   & 0                          & 9407                     \\
                     & 3                         &  t.l.                      & 1.38                       & 532                     &  & \textbf{2253.0}          & 0                          & 2221                    &  & 2472.1                   & 0                          & 7191                     \\
                     & 5                         &  t.l.                      & 1.25                       & 514                     &  & T.L.                     & 0                          & 5172                    &  & \textbf{3498.2}          & 0                          & 15334                    \\
                     & 7                         &  t.l.                     & 1.16                       & 496                     &  & 2742.4                   & 0                          & 4771                    &  & \textbf{1253.5}          & 0                          & 4224                     \\
                     & 10                        &  t.l.                      & 1.07                       & 491                     &  & 2170.8                   & 0                          & 2295                    &  & \textbf{1872.0}          & 0                          & 7258                     \\
                     & 15                        &  t.l.                      & 1.11                       & 506                     &  & 2846.1                   & 0                          & 5190                    &  & \textbf{910.5}           & 0                          & 3103                     \\
                     & 20                        &  t.l.                      & 0.72                       & 612                     &  & 2030.5                   & 0                          & 2419                    &  & \textbf{837.7}           & 0                          & 2880                     \\
\bottomrule
\end{tabular}
\end{table}

According to both tables,  {\fontfamily{qcr}\selectfont{CVX-D}} requires significantly longer computational time than other approaches. In particular, Table~\ref{Tab:Result200-1} shows that {\fontfamily{qcr}\selectfont{MICQP-D}} and {\fontfamily{qcr}\selectfont{MICQP-A}} are on average 3-5 times faster than  {\fontfamily{qcr}\selectfont{CVX-D}}. Furthermore, the instances under $\xi = 0.5$ are significantly harder than those under $\xi =0.1$. As shown in Table~\ref{Tab:Result200-05},  {\fontfamily{qcr}\selectfont{CVX-D}} can only solve 1 instance (out of 16) to optimality, whereas {\fontfamily{qcr}\selectfont{MICQP-A}} can solve all instances within the time limit. Therefore, we conclude that our proposed conic approaches are more efficient than the benchmark. In effect, as in Section~\ref{S:micqp}, it is rather straightforward to reformulate the model into a MICQP. We thus recommend MICQP-based approaches to solve the LLPTL.

We then proceed to compare the performance of {\fontfamily{qcr}\selectfont{MICQP-D}} and {\fontfamily{qcr}\selectfont{MICQP-A}}. From Table~\ref{Tab:Result200-1}-\ref{Tab:Result200-05}, {\fontfamily{qcr}\selectfont{MICQP-A}} performs better than {\fontfamily{qcr}\selectfont{MICQP-D}} in terms of the average computational time. However, we observe that when $\gamma$ is small, {\fontfamily{qcr}\selectfont{MICQP-D}} typically runs faster than {\fontfamily{qcr}\selectfont{MICQP-A}}, i.e., the computational time by {\fontfamily{qcr}\selectfont{MICQP-D}} is, in general, shorter than that by {\fontfamily{qcr}\selectfont{MICQP-A}} when $\gamma = 1,2,3$. This could be explained by the fact that {\fontfamily{qcr}\selectfont{MICQP-A}} uses the ADC constraint to model the dominance. Such a constraint aggregates a set of the DDCs into a single ADC. When $\gamma$ is small, the number of DDCs increases, and the resultant ADC needs to aggregate more constraints and thus leads to a weaker continuous relaxation according to Lemma~\ref{prop:DDC_ADC}. It turns out that despite having a substantially smaller constraint size, the advantage of {\fontfamily{qcr}\selectfont{MICQP-A}} in terms of the problem scale is overwhelmed by the disadvantages of its weaker relaxation and more required branching (see the column ``\#N"). However, when $\gamma$ is large (i.e., $\gamma \geq 5$), {\fontfamily{qcr}\selectfont{MICQP-A}} generally outperforms {\fontfamily{qcr}\selectfont{MICQP-D}} as we observe considerable speed-up by using the ADC in most of the reported instances.

\subsection{Result analysis on DS2}

Next, we investigate the performance of the MICQP approaches on the DS2 instances. Since these instances are significantly more challenging than previous ones, we set the optimality gap tolerance to 0.1\% and the time limit to 7200 seconds. Table~\ref{Tab:DS2} reports the results. For each problem instance, the computational time for the approach that performs the best is highlighted in boldface. For those instances where both approaches cannot terminate optimally within the time limit, the smaller $gp$ at termination is highlighted in boldface.

Table~\ref{Tab:DS2} suggests that {\fontfamily{qcr}\selectfont{MICQP-A}} is in general more efficient since it beats {\fontfamily{qcr}\selectfont{MICQP-D}} in most of the instances (27 out of 32). Meanwhile, both approaches show outstanding capabilities in solving these large-scale instances: Under $\xi = 1$, all instances are successfully solved  within the tolerance in 7200 seconds; Under $\xi = 0.5$, although there are unsolved instances, $gp$ is less than 0.5\% for {\fontfamily{qcr}\selectfont{MICQP-D}} and  0.3\% for {\fontfamily{qcr}\selectfont{MICQP-A}}. Both values are sufficiently small to guarantee the quality of the solutions. This experiment thus demonstrates that the proposed MICQP approaches are practically implementable.

\begin{table}[h]
\footnotesize
\centering
\caption{Computational results of the DS2 instances. Optimality Tolerance is 0.1\%. Time limit is 7200 seconds.}
\label{Tab:DS2}
\begin{tabular}{lrrrrrrrrrrrr} 
\toprule
\multirow{3}{*}{$f$} & \multirow{3}{*}{$\gamma$} & \multicolumn{5}{c}{$\xi=0.5$}                                                                                                        & ~ ~~                       & \multicolumn{5}{c}{$\xi=1.0$}                                                                                     \\ 
\cline{3-7}\cline{9-13}
                     &                           & \multicolumn{2}{c}{{\fontfamily{qcr}\selectfont{MICQP-D}}}                           & \multicolumn{1}{l}{} & \multicolumn{2}{c}{{\fontfamily{qcr}\selectfont{MICQP-A}}}                           & \multicolumn{1}{l}{~ ~ ~~} & \multicolumn{2}{c}{{\fontfamily{qcr}\selectfont{MICQP-D}}}                           &  & \multicolumn{2}{c}{{\fontfamily{qcr}\selectfont{MICQP-A}}}                            \\ 
\cline{3-4}\cline{6-7}\cline{9-10}\cline{12-13}
                     &                           & \multicolumn{1}{c}{t[s]} & \multicolumn{1}{c}{gp[\%]} &                      & \multicolumn{1}{c}{t[s]} & \multicolumn{1}{c}{gp[\%]} & ~ ~                        & \multicolumn{1}{c}{t[s]} & \multicolumn{1}{c}{gp[\%]} &  & \multicolumn{1}{c}{t[s]} & \multicolumn{1}{c}{gp[\%]}  \\ 
\hline

500                  & 1                         & \textbf{479.5}           & $<0.1$                          &                      & 1021.3                   & $<0.1$                          &                                                                                            & \textbf{205.2}           & $<0.1$                          &  & 582.2                    & $<0.1$                           \\
                     & 2                         & 1349.8                   &  $<0.1$                         &                      & \textbf{948.5}           & $<0.1$                          &                                                                                            & 481.9                    & $<0.1$                          &  & \textbf{471.7}           & $<0.1$                           \\
                     & 3                         & 1080.6                   & $<0.1$                          &                      & \textbf{1067.7}          & $<0.1$                          &                                                                                            & \textbf{282.9}           & $<0.1$                         &  & 522.1                    & $<0.1$                           \\
                     & 5                         & 1448.3                   & $<0.1$                          &                      & \textbf{958.8}           & $<0.1$                         &                                                                                            & 615.9                    & $<0.1$                          &  & \textbf{375.7}           & $<0.1$                           \\
                     & 7                         & 1703.6                   & $<0.1$                          &                      & \textbf{1151.4}          & $<0.1$                          &                                                                                            & 544.5                    & $<0.1$                          &  & \textbf{363.6}           & $<0.1$                           \\
                     & 10                        & 1773.8                   & $<0.1$                          &                      & \textbf{949.8}           & $<0.1$                          &                                                                                            & 631.9                    & $<0.1$                          &  & \textbf{383.5}           & $<0.1$                           \\
                     & 15                        & 1670.1                   & $<0.1$                          &                      & \textbf{482.7}           & $<0.1$                          &                                                                                            & 602.7                    & $<0.1$                          &  & \textbf{309.5}           & $<0.1$                           \\
                     & 20                        & 1781.5                   & $<0.1$                          &                      & \textbf{551.8}           & $<0.1$                          &                                                                                            & 543.0                    & $<0.1$                          &  & \textbf{192.4}           & $<0.1$                           \\
                     &                           &                          &                            &                      &                          &                            &                                                                                            &                          &                            &  &                          &                             \\
1000                 & 1                         & \textbf{1226.1}          & $<0.1$                          &                      & 1782.1                   & $<0.1$                          &                                                                                            & \textbf{848.0}           & $<0.1$                          &  & 1343.9                   & $<0.1$                           \\
                     & 2                         &  t.l.                      & 0.16                       &                      & \textbf{4219.2}          & $<0.1$                          &                                                                                            & 2955.8                   & $<0.1$                          &  & \textbf{1528.4}          & $<0.1$                           \\
                     & 3                         &  t.l.                     & 0.46                       &                      &  t.l.                     & \textbf{0.17}              &                                                                                            & 3918.4                   & $<0.1$                         &  & \textbf{2273.6}          & $<0.1$                           \\
                     & 5                         &  t.l.                      & 0.47                       &                      &  t.l.                      & \textbf{0.28}              &                                                                                            & 6044.9                   & $<0.1$                          &  & \textbf{1973.8}          & $<0.1$                           \\
                     & 7                         &  t.l.                      & 0.46                       &                      &  t.l.                      & \textbf{0.22}              &                                                                                            & 4642.7                   & $<0.1$                          &  & \textbf{4175.7}          & $<0.1$                           \\
                     & 10                        &  t.l.                      & 0.48                       &                      &  t.l.                      & \textbf{0.19}              &                                                                                            & 4787.9                   & $<0.1$                          &  & \textbf{2935.1}          & $<0.1$                           \\
                     & 15                        &  t.l.                      & 0.36                       &                      &  t.l.                      & \textbf{0.23}              &                                                                                            & 5274.9                   & $<0.1$                          &  & \textbf{2557.6}          & $<0.1$                           \\
                     & 20                        &  t.l.                      & 0.38                       &                      &  t.l.                      & \textbf{0.21}              &                                                                                            & 4677.5                   & $<0.1$                          &  & \textbf{2095.1}          & $<0.1$                           \\
\bottomrule
\end{tabular}
\end{table}

\section{Numerical study and implication}
\label{S:numerical}
In this section, we conduct sensitivity analysis to assess the impact of problem parameters on the solution structure and draw implications. The problem data used in our discussion is the DS1 dataset, and the attraction functions are defined according to (\ref{eqt:attr_locker}) and (\ref{eqt:attr_outside}), where parameter $\alpha \geq 0 $ reflects customers’ sensitivity to the distance and $\xi > 0$ controls the magnitude of the attraction of the outside option (i.e., $a_{i0}, \forall i \in I$).

We assume that the attraction of a locker to a customer zone only depends on the distance. In practice, other factors may also be important. \cite{lyu2019last} conducted a study that calibrates the utility. They combined two sources of delivery datasets, which contain the traditional home delivery records (including customer address, delivery time, delivery status) and the delivery to lockers records (including customer address, locker location, delivery time, customer pickup time, delivery status). Regression models were then developed to fit the data and derived the utility/attraction function. Unfortunately, this paper does not have access to the delivery data and thus, we link the attraction to the distance. Nevertheless, we should note that our locker location models are ready to generalize to more complicated utility or attraction functions.

\subsection{Comparing TLM with BNL and MNL}


To start with, we investigate the impact of $\gamma$ on the solution structure. We fix $\alpha=1$ and vary $f$, $\xi$, and $\gamma$. Notice that when $\gamma = 0$, the TLM becomes the BNL since the resultant dominance relation implies that customers will only consider the most attractive locker facility and the outside option. When $\gamma=\infty$, the TLM is exactly the MNL since all open lockers will in the nondominated set and thus have positive probability of being selected. Therefore, by varying the value of $\gamma$, we can compare the solutions under different TLM settings with the solutions under both the BNL and the MNL. From here onwards, we use TLM-$\gamma$ to denote the LLPTL model under a specific $\gamma$ value.

Table~\ref{Tab:compare_DCM} shows the numerical results under different choice models. Here $R$ and $\#F$ stand for the estimated revenue and the number of open locker facilities.  According to the table, for all combinations of $(f,\xi)$, both the profit and the number of open locker facilities increase with $\gamma$ (i.e., Profit and $\#F$ increase from the ``BNL" column to the ``MNL" column). Intuitively, the number of dominated lockers for a customer zone decreases with $\gamma$. With a higher $\gamma$ value, customers are more receptive to the locker service and are willing to use more lockers as their delivery options. Therefore, by opening more lockers, the company expects to capture more demand (i.e., yield a higher revenue) and eventually leads to a higher profit compared to the cases with lower $\gamma$ values.

\eat{

\begin{table}[h]
\small
\centering
\caption{Solution structures under different choice models.}
\label{Tab:compare_DCM}
\begin{tabular}{lllllllll} 
\toprule
$(f,\xi)$                   &          & \textbf{BNL} & \textbf{TLM-1} & \textbf{TLM-2} & \textbf{TLM-3} & \textbf{TLM-5} & \textbf{TLM-10} & \textbf{MNL}  \\ 
\hline
\multirow{4}{*}{(500,0.5)}  & Profit   & 28737        & 31389         & 32613         & 33138         & 34010         & 34804          & 36380         \\
                            & \#F      & 32           & 37            & 39            & 40            & 40            & 40             & 42            \\
                            & R        & 44737        & 49889         & 52113         & 53138         & 54010         & 54804          & 57380         \\
                            & $\Delta[\%]$ & 28.26        & 15.02         & 10.11         & 7.98          & 6.24          & 4.70           & -             \\
                            &          &              &               &               &               &               &                &               \\
\multirow{4}{*}{(500,1.0)}  & Profit   & 15080        & 16873         & 17842         & 18303         & 18933         & 19601          & 21037         \\
                            & \#F      & 28           & 32            & 34            & 36            & 37            & 38             & 41            \\
                            & R        & 29080        & 32873         & 34842         & 36303         & 37433         & 38601          & 41537         \\
                            & $\Delta[\%]$ & 42.84        & 26.36         & 19.22         & 14.42         & 10.96         & 7.61           & -             \\
                            &          &              &               &               &               &               &                &               \\

\multirow{4}{*}{(1000,0.5)} & Profit   & 16045        & 17001         & 17517         & 17805         & 18317        & 18994          & 20004         \\
                            & \#F      & 19           & 20            & 23            & 23            & 23            & 24             & 25            \\
                            & R        & 35045        & 37001         & 40517         & 40805         & 41994         & 42994          & 45004         \\
                            & $\Delta[\%]$ & 28.42        & 21.63         & 11.07         & 10.29         & 8.92          & 4.68           & -             \\
                            &          &              &               &               &               &               &                &               \\

\multirow{4}{*}{(1000,1.0)} & Profit   & 5795         & 6202          & 6328          & 6411          & 6635          & 6901           & 7442          \\
                            & \#F      & 11           & 13            & 14            & 14            & 15            & 16             & 18            \\
                            & R        & 16795        & 19202         & 20328         & 20411         & 21635         & 22901          & 25442         \\
                            & $\Delta[\%]$  & 51.49        & 32.50         & 25.16         & 24.65         & 17.60         & 11.10          & -             \\
\bottomrule
\end{tabular}
\end{table}
}

\begin{table}
\footnotesize
\centering
\caption{Solution structures under different choice models.}
\label{Tab:compare_DCM}
\begin{tabular}{llllllll} 
\toprule
$(f,\xi)$                   &              & \textbf{BNL} & \textbf{TLM-1} & \textbf{TLM-2} & \textbf{TLM-3} & \textbf{TLM-5} & \textbf{MNL}  \\ 
\hline
\multirow{4}{*}{(500,0.5)}  & Profit       & 28737        & 31389          & 32613          & 33138          & 34010          & 36380         \\
                            & \#F          & 32           & 37             & 39             & 40             & 40             & 42            \\
                            & R            & 44737        & 49889          & 52113          & 53138          & 54010          & 57380         \\
                            & $\Delta[\%]$ & 28.26        & 15.02          & 10.11          & 7.98           & 6.24           & -             \\
                            &              &              &                &                &                &                &               \\
\multirow{4}{*}{(500,1.0)}  & Profit       & 15080        & 16873          & 17842          & 18303          & 18933          & 21037         \\
                            & \#F          & 28           & 32             & 34             & 36             & 37             & 41            \\
                            & R            & 29080        & 32873          & 34842          & 36303          & 37433          & 41537         \\
                            & $\Delta[\%]$ & 42.84        & 26.36          & 19.22          & 14.42          & 10.96          & -             \\
                            &              &              &                &                &                &                &               \\
\multirow{4}{*}{(1000,0.5)} & Profit       & 16045        & 17001          & 17517          & 17805          & 18317          & 20004         \\
                            & \#F          & 19           & 20             & 23             & 23             & 23             & 25            \\
                            & R            & 35045        & 37001          & 40517          & 40805          & 41994          & 45004         \\
                            & $\Delta[\%]$ & 28.42        & 21.63          & 11.07          & 10.29          & 8.92           & -             \\
                            &              &              &                &                &                &                &               \\
\multirow{4}{*}{(1000,1.0)} & Profit       & 5795         & 6202           & 6328           & 6411           & 6635           & 7442          \\
                            & \#F          & 11           & 13             & 14             & 14             & 15             & 18            \\
                            & R            & 16795        & 19202          & 20328          & 20411          & 21635          & 25442         \\
                            & $\Delta[\%]$ & 51.49        & 32.50          & 25.16          & 24.65          & 17.60          & -             \\
\bottomrule
\end{tabular}
\end{table}

To directly quantify the impact of using different choice models, we use $\Delta$ to denote the relative difference in the estimated revenue. It is computed by 
\begin{align}
\Delta = \frac{\text{R(MNL)} - \text{R(TLM-$\gamma$)}}{\text{R(TLM-$\gamma$)}} \times 100\%
\end{align}
where R(MNL) and  R(TLM-$\gamma$) are the estimated revenue under the MNL and the TLM-$\gamma$ (note that the TLM-$0$ is the BNL). $\Delta$ thus reflects to what extent the MNL results in higher estimated revenue compared to other models.

From Table~\ref{Tab:compare_DCM}, the MNL reports significantly higher estimated revenue compared to the BNL: the corresponding $\Delta$ value exceeds 28\% and blows up to more than 51\% when $f=1000$ and $\xi=1.0$. Such a result is owing to the different modeling assumptions. Customers under the MNL maintain all open lockers in their consideration set, and thus the attraction of the locker service depends on all lockers (see the expression in Equation~(\ref{eqt:p_tlm})). By contrast, the BNL assumes that customers only consider the open locker with the highest attraction. Consequently, the attraction of the service under the BNL could be substantially lower, which naturally leads to a lower estimated revenue because more customers are predicted to select the outside option.

\begin{figure*}[h]
~~~~~ \psfig{figure=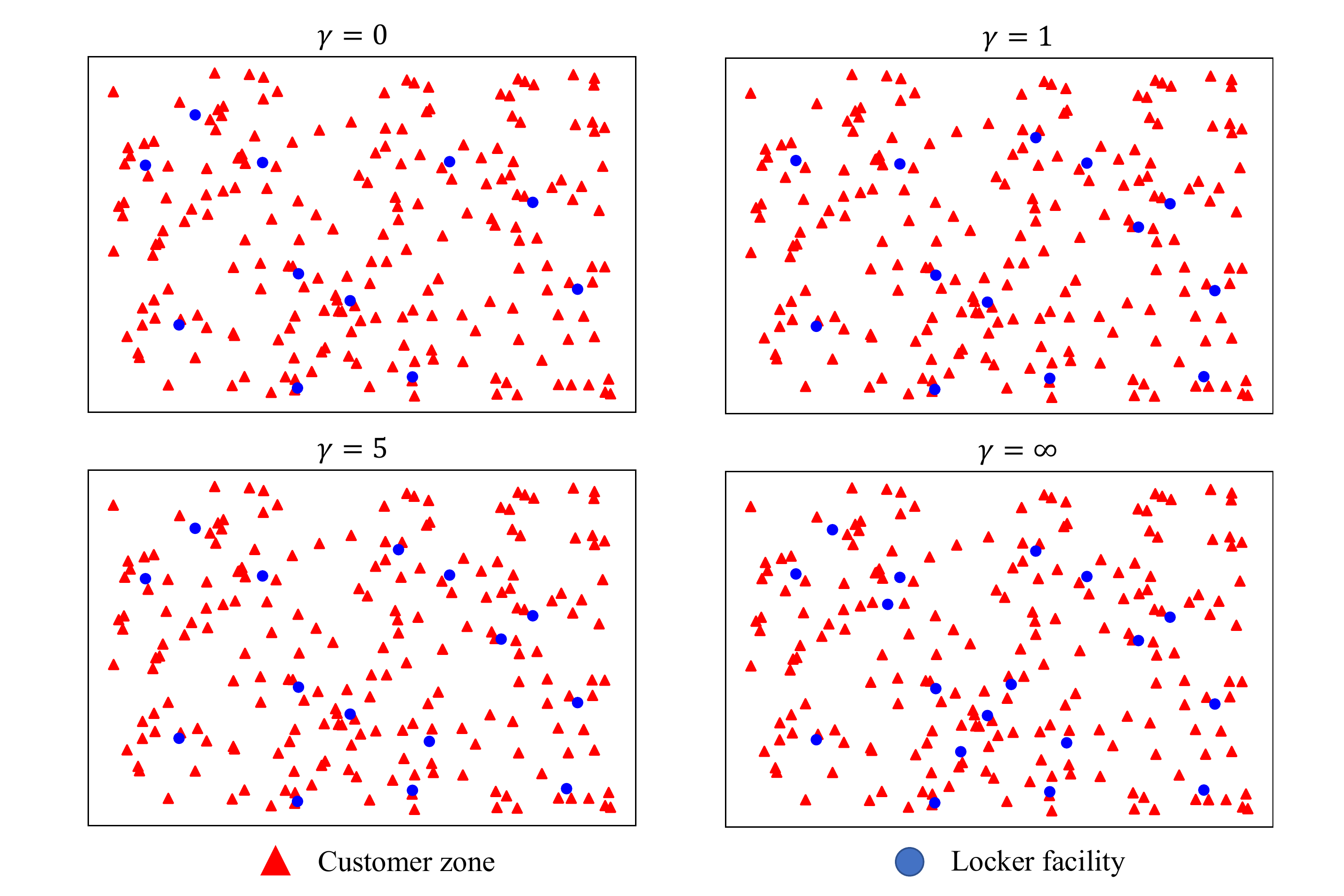,width=150mm}
\caption{Geographic locations of customer zones and locker facilities. $(f,\xi) = (1000,1.0)$.} \label{fig:geo_plot}
\end{figure*}

In effect, the assumptions on the BNL and the MNL are rather strong because  the models assume that customers visit either only one facility or all facilities. Both scenarios are essentially special cases of a more generalized framework. If customer's choice behaviors follow the TLM, then using the BNL and the MNL could, respectively, \textit{overestimate} and \textit{underestimate} the revenue. In both cases, the resultant location decisions are clearly suboptimal, and their difference is so remarkable that it is hard to justify the quality of the revenue estimation. For example, when $f=1000$ and $\xi=1.0$, the number of open facilities under the MNL is 18, which is 7 more than that under the BNL (for geographical location of facilities, see Figure~\ref{fig:geo_plot}), and thus incurs additional facilities cost of 7000. Considering that the estimated profit under the MNL is only 7442, such an additional investment is relatively extensive. In a nutshell, the MNL in the locker location problem would lead to an \textit{aggressive} decision that suggests to opening a larger number of facilities, whereas the decision by the BNL would be  \textit{conservative}. 

Due to the strong assumptions on the BNL and the MNL, it is essential to develop an intermediate and generalized approach to conduct a more accurate revenue estimation. This motivates the use of TLM.  As in Table~\ref{Tab:compare_DCM}, when $\gamma$ increases, the TLM shift away from the BNL to the MNL, and we observe a decreasing $\Delta$ and an increasing $\#F$. The increased number of open facilities can also be directly seen from Figure~\ref{fig:geo_plot}. In general, by setting an appropriate value of $\gamma$, we can effectively balance the trade-off between aggressive and conservative decisions and avoid over-pessimistic and over-optimistic estimation of the revenue.  

\subsection{Profit loss due to aggressive decision}

In the previous experiment, we demonstrate that using the MNL will overestimate the revenue. As a result, the model will suggest opening too many facilities. This experiment will investigate the consequence, i.e., the profit loss, due to such an aggressive decision.

To facilitate our discussion, we introduce the following notations. Let $x_{\infty}$ be the location decision when we employ the MNL. The \textit{Actual Profit} under the TLM-$\gamma$ is evaluated by holding the location decision at $x_{\infty}$ and optimizing the remaining problem of MICQP. Therefore,  \textit{Actual Profit} stands for the profit when the decision is made based on the MNL but the actual choice model is TLM-$\gamma$. Moreover, the  \textit{Optimal Profit} is the optimal objective function value when we use the correct choice model. Then, we define \textit{RelLoss}[\%] as
\begin{align}
RelLoss = \frac{Optimal~Profit-Actual~Profit}{Optimal~Profit} \times 100\%
\end{align}
which specifies the relatively profit loss due to the use of the MNL when the choice model does not follow the MNL.

Table~\ref{Tab:profit_loss} reports the results of the profit loss. We observe that $RelLoss$ reduces when $\gamma$ increases. This is reasonable and intuitive since a larger $\gamma$ value implies that the choice model behaves more similarly to the MNL, and thus the location decision made based on the MNL is closer to that based on the TLM-$\gamma$. However, When $\gamma$ is small, $RelLoss$ is significant. In particular, when $\gamma=0$ (i.e., the BNL), $RelLoss$ can below up to more than 40\%, and  even when $\gamma=2$, $RelLoss$ is still on a relatively high value for some combinations of $(f,\xi,\alpha)$. These results illustrate the importance of considering the zero-choice probability in the locker location problem, suggesting that an aggressive decision due to the use of the MNL will incur an unnecessarily high facility cost that cannot be compensated by the additional revenue.

\begin{table}[h]
\centering
\footnotesize
\caption{Profit loss due to the use of the MNL.}
\label{Tab:profit_loss}
\setlength{\tabcolsep}{0.8mm}{
\begin{tabular}{lllllllllllll} 
\toprule
\multirow{2}{*}{$(f,\xi)$}  &                & \multicolumn{5}{c}{$\alpha=0.8$}                                                 & ~~ & \multicolumn{5}{c}{$\alpha=1.0$}                                                                                                   \\ 
\cline{3-7}\cline{9-13}
                            &                & \textbf{BNL} & \textbf{TLM-1} & \textbf{TLM-2} & \textbf{TLM-3} & \textbf{TLM-5} &  & \textbf{\textbf{BNL}} & \textbf{\textbf{TLM-1}} & \textbf{\textbf{TLM-2}} & \textbf{\textbf{TLM-3}} & \textbf{\textbf{TLM-5}}  \\ 

\hline
\multirow{3}{*}{(500,0.5)}  & Optimal Profit & 32656        & 36599          & 37935          & 38903          & 39931          &  & 28737                 & 31389                    & 32613                    & 33138                    & 34010                     \\
                            & Actual Profit  & 30562        & 35761          & 37619          & 38714          & 39867          &  & 26725                 & 30314                    & 32222                    & 32967                    & 33985                     \\
                            & RelLoss$[\%]$ & 6.41         & 2.29           & 0.83           & 0.49           & 0.16           &  & 7.00                  & 3.42                     & 1.20                     & 0.52                     & 0.07                      \\
                            &                &              &                &                &                &                &  &                       &                          &                          &                          &                           \\
\multirow{3}{*}{(500,1.0)}  & Optimal Profit & 17824        & 20923          & 21998          & 22863          & 24030          &  & 15080                 & 16873                    & 17842                    & 18303                    & 18933                     \\
                            & Actual Profit  & 13633        & 19561          & 21373          & 22407          & 23770          &  & 12638                 & 15992                    & 17439                    & 18035                    & 18805                     \\
                            & RelLoss$[\%]$ & 23.51        & 6.51           & 2.84           & 1.99           & 1.08           &  & 16.19                 & 5.22                     & 2.26                     & 1.46                     & 0.68                      \\
                            &                &              &                &                &                &                &  &                       &                          &                          &                          &                           \\
\multirow{3}{*}{(500,1.5)}  & Optimal Profit & 10912        & 13023          & 13864          & 14501          & 15418          &  & 8760                  & 9834                     & 10579                    & 10940                    & 11459                     \\
                            & Actual Profit  & 6560         & 11620          & 13142          & 13965          & 15128          &  & 6670                  & 9260                     & 10311                    & 10834                    & 11366                     \\
                            & RelLoss$[\%]$ & 39.88        & 10.77          & 5.21           & 3.70           & 1.88           &  & 23.86                 & 5.84                     & 2.53                     & 0.97                     & 0.81                      \\
                            &                &              &                &                &                &                &  &                       &                          &                          &                          &                           \\
\multirow{3}{*}{(1000,0.5)} & Optimal Profit & 20592        & 22549          & 23097          & 23622          &    24506            &  & 16045                 & 17001                    & 17517                    & 17805                    & 18317                     \\
                            & Actual Profit  & 18539        & 21677          & 22384          & 23175          & 24224          &  & 14865                 & 16551                    & 17416                    & 17690                    & 18249                     \\
                            & RelLoss$[\%]$ & 9.97         & 3.87           & 3.09           & 1.89           &    1.15            &  & 7.35                  & 2.65                     & 0.58                     & 0.65                     & 0.37                      \\
                            &                &              &                &                &                &                &  &                       &                          &                          &                          &                           \\
\multirow{3}{*}{(1000,1.0)} & Optimal Profit & 8559         & 9436           & 9670           & 10025          & 10557          &  & 5795                  & 6202                     & 6328                     & 6411                     & 6635                      \\
                            & Actual Profit  & 6556         & 8302           & 8641           & 9290           & 10099          &  & 4851                  & 5426                     & 5847                     & 5971                     & 6425                      \\
                            & RelLoss$[\%]$ & 23.40        & 12.02          & 10.64          & 7.33           & 4.34           &  & 16.29                 & 12.51                    & 7.60                     & 6.86                     & 3.17                      \\
                            &                &              &                &                &                &                &  &                       &                          &                          &                          &                           \\
\multirow{3}{*}{(1000,1.5)} & Optimal Profit & 3969         & 4282           & 4403           & 4487           & 4637           &  & 2426                  & 2472                     & 2501                     & 2517                     & 2533                      \\
                            & Actual Profit  & 2286         & 3354           & 3586           & 3903           & 4360           &  & 2198                  & 2309                     & 2361                     & 2393                     & 2447                      \\
                            & RelLoss$[\%]$ & 42.40        & 21.67          & 18.56          & 13.02          & 5.97           &  & 9.40                  & 6.59                     & 5.60                     & 4.93                     & 3.40                      \\
\bottomrule

\end{tabular}}
\end{table}

\eat{
\begin{table}[h]
\centering
\caption{Profit loss due to the use of the MNL.}
\label{Tab:compare_DCM}
\begin{tabular}{llllllll}
\toprule
$(f,\xi)$                   &                & \textbf{BNL} & \textbf{TLM-1} & \textbf{TLM-2} & \textbf{TLM-3} & \textbf{TLM-5} & \textbf{TLM-10}  \\ 
\hline
\multirow{3}{*}{(500,0.5)}  & Optimal Profit & 28737        & 31389          & 32613          & 33138          & 34010          & 34804            \\
                            & Actual Profit  & 26725        & 30314          & 32222          & 32967          & 33985          & 34785            \\
                            & Rel-Loss$[\%]$ & 7.00         & 3.42           & 1.20           & 0.52           & 0.07           & 0.05             \\
                            &                &              &                &                &                &                &                  \\
\multirow{3}{*}{(500,1.0)}  & Optimal Profit & 15080        & 16873          & 17842          & 18303          & 18933          & 19601            \\
                            & Actual Profit  & 12638        & 15992          & 17439          & 18035          & 18805          & 19521            \\
                            & Rel-Loss$[\%]$ & 16.19        & 5.22           & 2.26           & 1.46           & 0.68           & 0.41             \\
                            &                &              &                &                &                &                &                  \\
\multirow{3}{*}{(500,1.5)}  & Optimal Profit & 8760         & 9834           & 10579          & 10940          & 11459          & 11983            \\
                            & Actual Profit  & 6670         & 9260           & 10311          & 10834          & 11366          & 11903            \\
                            & Rel-Loss$[\%]$ & 23.86        & 5.84           & 2.53           & 0.97           & 0.81           & 0.67             \\
                            &                &              &                &                &                &                &                  \\
\multirow{3}{*}{(1000,0.5)} & Optimal Profit & 16045        & 17001          & 17517          & 17805          & 18317          & 18994            \\
                            & Actual Profit  & 14865        & 16551          & 17416          & 17690          & 18249          & 18971            \\
                            & Rel-Loss$[\%]$ & 7.35         & 2.65           & 0.58           & 0.65           & 0.37           & 0.12             \\
                            &                &              &                &                &                &                &                  \\                            
\multirow{3}{*}{(1000,1.0)} & Optimal Profit & 5795         & 6202           & 6328           & 6411           & 6635           & 6901             \\
                            & Actual Profit  & 4851         & 5426           & 5847           & 5971           & 6425           & 6830             \\
                            & Rel-Loss$[\%]$ & 16.29        & 12.51          & 7.60           & 6.86           & 3.17           & 1.03             \\
                            &                &              &                &                &                &                &                  \\                            
\multirow{3}{*}{(1000,1.5)} & Optimal Profit & 2426         & 2472           & 2501           & 2517           & 2533           & 2583             \\
                            & Actual Profit  & 2198         & 2309           & 2361           & 2393           & 2447           & 2571             \\
                            & Rel-Loss$[\%]$ & 9.40         & 6.59           & 5.60           & 4.93           & 3.40           & 0.46  \\ 
\bottomrule                                     
\end{tabular}
\end{table}
}

\subsection{Impact of $\alpha$}

The next experiment studies how $\alpha$, i.e., customers’ sensitivity to the distance, affects the solution. We fix $F=500$ and $\xi = 1$ and vary the value of $\alpha$. 

\begin{figure}[h]
\begin{center}
\subfigure[Profit versus $\alpha$.]{
     \psfig{figure=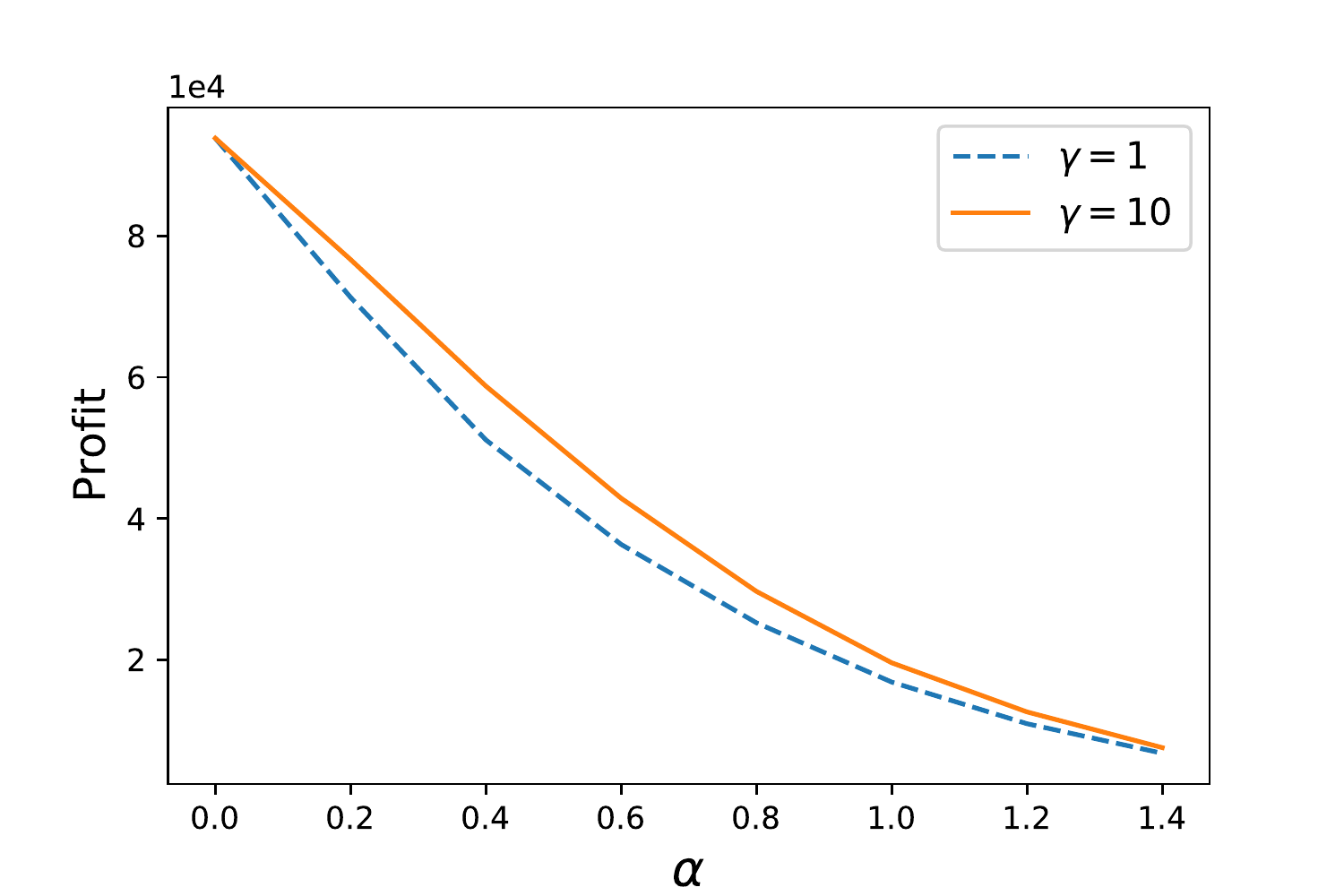,width=78mm} \label{fig:case_profit_alpha}
     }
\subfigure[\#F versus $\alpha$.]{
     \psfig{figure=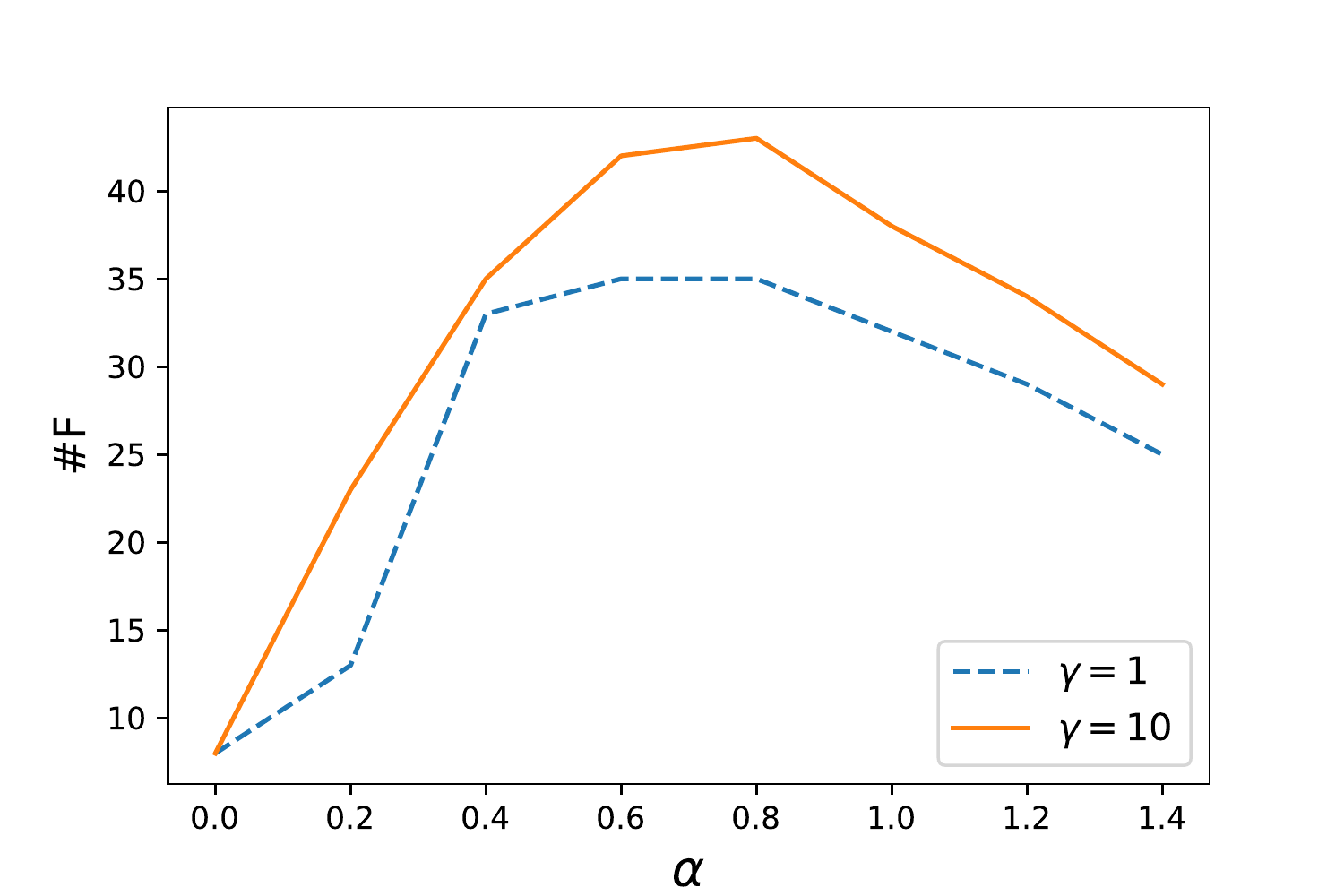,width=78mm} \label{fig:case_F_alpha}
     }     
\end{center}
\caption{The impact of $\alpha$ on the profit and the number of open locker facilities.} \label{alpha}
\end{figure}

Figure~\ref{alpha} depicts the results under $\gamma=1$ and $\gamma=10$. We observe that the profit decreases with $\alpha$. Intuitively, under a large $\alpha$ value, the attraction of the locker facilities is  low, and thus the attraction of the outside option becomes relatively high, implying that customers are more willing to use the outside option. By contrast, if $\alpha$ is small, then customers are not sensitive to the distance and perceive high utility of the locker service. In this case, more customers will select the parcel locker as their delivery mode, leading to a high profit of the company.  

Interestingly, the number of open facilities $\#F$ first increases with $\alpha$ to a certain maximum and then decreases. When $\alpha$ is at a very small value, the attraction of each locker is relatively high. The company can attract a large proportion of customers by opening only a small number of facilities. Now, when $\alpha$ increases (from 0 to around 0.8), the locker service becomes less attractive but stays competitive compared to the outside option. Then the company needs to open more facilities so as to attract more customers and improve the revenue. However, if customers are very sensitive to the distance (i.e., $\alpha > 1$), then the attraction of the locker service is rather low, and customers may view the locker service as a less competitive option and prefer not to use it. For the company, the increased revenue by opening additional facilities in such a situation will not compensate for the cost; therefore, it is recommended to open fewer facilities to save the investment when $\alpha$ keeps growing.

\eat{
\begin{table}
\centering
\caption{$|I|=200, |J|=100$, $\xi=0.5$}
\label{Tab:Result200-1}
\begin{tabular}{lllllllllll} 
\hline
\multirow{2}{*}{$(f,\gamma)$} &        & \multicolumn{9}{c}{$\xi$}                                                                                                              \\ 

                              &        & \textbf{0.05} & \textbf{0.1} & \textbf{0.3} & \textbf{0.5} & \textbf{0.7} & \textbf{1.0} & \textbf{1.3} & \textbf{1.5} & \textbf{2.0}  \\ 
\hline
\multirow{3}{*}{(500,1)}      & Profit &               & 63166        & 42421        & 31389        & 24041        & 16873        & 12150        & 9834         & 5926          \\
                              & R      &               & 80166        & 60921        & 49889        & 42541        & 32873        & 27650        & 23334        & 15926         \\
                              & \#F    &               & 34           & 37           & 37           & 37           & 32           & 31           & 27           & 20            \\
                              &        &               &              &              &              &              &              &              &              &               \\
\multirow{3}{*}{(500,10)}     & Profit & 74470         & 65472        & 45759        & 34804        & 27318        & 19601        & 14535        & 11983        & 7418          \\
                              & R      & 88970         & 82472        & 65259        & 54804        & 47318        & 38601        & 32035        & 27983        & 21918         \\
                              & \#F    & 29            & 34           & 39           & 40           & 40           & 38           & 35           & 32           & 29            \\
                              &        &               &              &              &              &              &              &              &              &               \\
\multirow{3}{*}{(1000,10)}    & Profit &               & 51154        & 29574        & 18994        & 12581        & 6901         & 3813         & 2583         & 930           \\
                              & R      &               & 75154        & 57574        & 42994        & 33581        & 22901        & 14813        & 10583        & 5930          \\
                              & \#F    &               & 24           & 28           & 24           & 21           & 16           & 11           & 8            & 5             \\
\hline
\end{tabular}
\end{table}
}

\subsection{Impact of $\xi$}

Finally, we conduct sensitivity analysis on $\xi$, i.e., the magnitude of the attraction of the outside option. We fix $\alpha=1$ and vary the value of $\xi$.

Table~\ref{Tab:xi} reports the results under 3 combinations of $(f,\gamma)$. One can observe that  increasing $\xi$ will decrease the profit. The reason is obvious: As the outside option becomes more attractive, more customers will select the outside option instead of the locker service. This naturally leads to a reduced revenue for the company.  It is worth noting that $\#F$ again shows a pattern where, as $\xi$ increases, the number of open facilities climbs up first and then declines.  Since the value of $\xi$ controls the attraction of the outside option and does not affect the attraction of the locker service, its impact on $\#F$ will be similar to that of $\alpha$, i.e., increasing either $\xi$ or $\alpha$ will make the outside option relatively more attractive (compared to the locker service). Therefore, the change of $\#F$ with respect to $\xi$ can be explained using the same observation in the last subsection.

\begin{table}[h]
\footnotesize
\centering
\caption{Impact of $\xi$ on profit, revenue, and the number of open facilities.}
\label{Tab:xi}
\begin{tabular}{llllllllllll} 
\toprule
\multirow{2}{*}{$\xi$} & \multicolumn{3}{c}{$f=500,\gamma=1$~~}                                       &  & \multicolumn{3}{c}{$f=500,\gamma=10$}                                        &  & \multicolumn{3}{c}{$f=1000,\gamma=10$}                                        \\ 
\cline{2-4}\cline{6-8}\cline{10-12}
                       & \multicolumn{1}{c}{Profit} & \multicolumn{1}{c}{R} & \multicolumn{1}{c}{\#F} &  & \multicolumn{1}{c}{Profit} & \multicolumn{1}{c}{R} & \multicolumn{1}{c}{\#F} &  & \multicolumn{1}{c}{Profit} & \multicolumn{1}{c}{R} & \multicolumn{1}{c}{\#F}  \\ 
\hline
0.05                   & 72868                      & 86868                 &  27                      &  & 74470                      & 88970                 & 29                      &  & 61988                           & 82988             & 21                          \\
0.1                    & 63166                      & 80166                 & 34                      &  & 65472                      & 82472                 & 34                      &  & 51154                      & 75154                 & 24                       \\
0.3                    & 42421                      & 60921                 & 37                      &  & 45759                      & 65259                 & 39                      &  & 29574                      & 57574                 & 28                       \\
0.5                    & 31389                      & 49889                 & 37                      &  & 34804                      & 54804                 & 40                      &  & 18994                      & 42994                 & 24                       \\
0.7                    & 24041                      & 42541                 & 37                      &  & 27318                      & 47318                 & 40                      &  & 12581                      & 33581                 & 21                       \\
1.0                    & 16873                      & 32873                 & 32                      &  & 19601                      & 38601                 & 38                      &  & 6901                       & 22901                 & 16                       \\
1.3                    & 12150                      & 27650                 & 31                      &  & 14535                      & 32035                 & 35                      &  & 3813                       & 14813                 & 11                       \\
1.5                    & 9834                       & 23334                 & 27                      &  & 11983                      & 27983                 & 32                      &  & 2583                       & 10583                 & 8                        \\
\bottomrule
\end{tabular}
\end{table}

To summarize, under different choice model settings ($\gamma$) and different ``relative" attraction of the locker service ($\alpha$ and $\xi$), the estimated revenue and the location decisions could be quite different. Therefore, it is important to carefully calibrate the parameters in the choice model (perhaps, through well-designed and structured questionnaires or statistical estimation from real big data) so that the company can better leverage the model to make decent decisions.

\section{Conclusion}
\label{S:conclusion}

This paper studied the location problem of a parcel locker system, taking customers' choices and zero-probability choice scenarios into considerations. More specifically, we employed the TLM to predict the revenue and proposed a combinatorial optimization problem to maximize the profit.  We then derived two equivalent mathematical formulations and proposed a unified solution approach based on mixed-integer conic quadratic programming.  Extensive numerical studies showed that the proposed methodology is practically implementable.  Finally, this paper investigated the impacts of parameters on the solution structure using sensitivity analysis. The results highlighted the necessity of considering the limited choice behavior of the customer and the importance of calibrating the customer choice model in the locker location problem.

There are limitations and future directions. We did not consider the locker capacity. In practice, capacity can be an important factor influencing customer choices. For example, if a locker has no remaining capacity, customers may select a nearby locker or opt for other delivery modes. However, with limited capacity, the TLM will fail because the probability distribution over the nondominated lockers will not follow the multinominal logit model anymore. In fact, when capacity comes into play, nearly all existing discrete choice models are not applicable. Indeed, It is challenging research problem to incorporate capacity constraints into the locker location problem with probabilistic choice behaviors.  Perhaps, one possible modeling approach to address this problem could be the entropy maximization rule; however, this will lead to a completely new study that positions itself outside the scope of this paper. Therefore, we would like to leave it for future research. Secondly, we proposed path-based inequalities to strengthen the ADC. However, we only derive the longest path and add one constraint for each customer zone. The representation of dominance relations as the directed acyclic graph could be further investigated to improve the computational efficiency. In particular, it is meaningful to explore how to leverage graph theory to derive strong inequalities and develop more efficient solution approaches for large-scale problems. Finally, we believe that an in-depth investigation on the factors that affect customer's choices (i.e., factors that define the utility) and the careful selection of choice models (including choice models such as mixed-logit and Markov chain choice model) will be useful when making location decisions, and these can also be an interesting research direction.

\section*{Acknowledgements}

This research has been made possible by the funding support from the Singapore Maritime Institute. The authors have also been supported by the Singapore A$*$STAR IAF-PP fund (Grant No. A1895a0033) under the project ``Digital Twin for Next Generation Warehouse". 


\bibliographystyle{apalike}
\bibliography{ref}

\end{document}